# Hazard Estimation under Generalized Censoring

Alberto Carabarin Aguirre*
albertocarabarin@yahoo.com.mx

B. Gail Ivanoff†
givanoff@uottawa.ca

Department of Mathematics and Statistics
University of Ottawa
585 King Edward Avenue
Ottawa, ON Canada K1N 6N5

November 1, 2018

*Research supported by Consejo Nacional de Ciencia y Tecnología, Mexico

†Research supported by a grant from the Natural Sciences and Engineering Research Council of Canada.

## Abstract

This paper focuses on the problem of the estimation of the cumulative hazard function of a distribution on a general complete separable metric space when the data points are subject to censoring by an arbitrary adapted random set. A problem involving observability of the estimator proposed in [8] and [9] is resolved and a functional central limit theorem is proven for the revised estimator. Several examples and applications are discussed, and the validity of bootstrap methods is established in each case.

# 1 Introduction

In two recent papers, [8] and [9], the problem of survival analysis on a general complete separable metric space was approached from the point of view of set-indexed martingales (cf. [6]). The goal of these papers was to build a mathematical structure that would handle data on general spaces (including, but not limited to $d$-dimensional Euclidean space) subject to very general types of censoring mechanisms. In particular, the cumulative hazard function was defined on a class of sets, and a corresponding Nelson-Aalen-type estimator was proposed.

In [8], the censoring mechanism in a set-indexed framework was defined by introducing the concept of a stopping set, a random set with a particular kind of measurability. The stopping sets acted as windows, since we could observe the data points only through them, and they permit the consideration of more sophisticated censoring schemes than the usual multivariate generalizations of right-censoring; however, the stopping sets were still restricted to certain shapes.

The theory of stopping sets was expanded to more general random sets called anti-clouds in [9]. The name is inspired by the example of aerial photography, where pictures are taken from high in the air and clouds may interfere with the observation of whatever the subject of interest is. Anti-clouds, then, are the complement of the clouds: the regions where we can actually observe the data points. The difference between the stopping sets and the anti-clouds is, quite simply, that anti-clouds have virtually no restrictions in the shape they can take; obviously, this makes for much more realistic censoring schemes.

The measurability requirement imposed on these random sets allow one to apply the theory of set-indexed martingales as developed in [6] to produce a set-indexed Nelson-Aalen estimator for the cumulative hazard on a general complete separable metric space in the presence of generalized censoring. In the case of censoring by stopping sets, consistency and asymptotic unbiasedness of the estimator, as well as asymptotic normality of its finite-dimensional distributions, are proven in [8]. The Nelson-Aalen estimator in the presence of censoring by clouds presented more of a challenge, due to the general geometric nature of the clouds; nevertheless, the estimator was shown to be consistent and asymptotically unbiased in [9].

Unfortunately, there were still some gaps in the general theory. Aside from the lack of a functional central limit theorem for the general estimator of [9], there remained a critical problem limiting the applicability of the estimator. As will be seen subsequently, a process needed to construct the Nelson-Aalen estimator may not be observable under some common data structures. In practice, this could render the estimator unusable much of the time.

Our first goal is to address the observability problem and produce a working estimator for the cumulative hazard, while still preserving the general censoring model of [9]. It happens that the root of the observability problem lies in the measurability requirement of the anti-clouds. By asking for a slightly weaker measurability condition, it is possible for us to achieve our objective. We call the new kind of random set that came from this modification a $*$-anti-cloud. If the complements of these sets, the $*$-clouds, act as the censoring in the survival model, then we are always able to construct a Nelson-Aalen estimate that can be observed.

Our second goal is to prove a functional central limit theorem for the Nelson-Aalen estimator under generalized censoring on Euclidean space. In this case, we can also establish the validity of bootstrap procedures.

It should be pointed out that various estimators for the cumulative hazard have been proposed for bivariate data under right censoring. However, for the most part they have been designed en route to producing a Kaplan-Meier-type estimator for the survival function. Unfortunately, in higher dimensions the relationship between the survival and the hazard functions is not as straightforward as it is in one dimension. In two dimensions, for example, the survival function is determined by the hazard and both marginal distributions. For a good exposition of the representation of a bivariate survival function $S$ in terms of the hazard and the marginals, we refer the reader to [5] or [10]. Usually, one-dimensional Kaplan-Meier estimates are used for the marginal survival functions; what will change is the way the hazard is estimated. As Kalbfleisch and Prentice explain in [10], the simplest of this type of estimator is due to Bickel, who uses the ratio of the number of failure points that were uncensored in both components at $(t_1, t_2)$ to the number of pairs at risk at that point. The Dabrowska estimator and the Prentice-Cai estimator represent attempts at improving this approach, by considering not only the ratio of the number of double failure times to the at-risk set, but also the ratio of the number of failures on one coordinate when the other is still alive to the at-risk set. Lastly, Pons proposed a martingale-based estimator of the cumulative hazard of two right-censored survival times in [11], which she then used to produce a test for independence of the two survival times by comparing it to the product of the usual one-dimensional Nelson-Aalen estimates of both marginal cumulative hazards. However, her model required independence of the components of the survival time.

All of the preceding estimators depend strongly on the structure of Euclidean space and as well on the censoring mechanism, which is generally assumed to be right censoring of each component of the data point separately. The advantage of the martingale methodology used here is that it provides a unified and versatile approach to hazard estimation; in particular it is applicable to both Euclidean and non-Euclidean spaces, and may be applied to very general censoring mechanisms.

We will proceed as follows. In §2 we provide the mathematical framework for our model. We develop the notion of the *-anti-cloud and show that certain martingale properties are preserved under filtering by this more general censoring mechanism. In §3, the observability problem of the original Nelson-Aalen estimator is discussed. The results of §2 allow us redefine the Nelson-Aalen estimator of the cumulative hazard under censoring by *-clouds. We present several examples to illustrate how our new estimator can be used to circumvent the observability problem through a transition from the anti-cloud model to the $*$-anti-cloud structure. The price we pay is a possible loss of previously uncensored data that will need to be discarded. In §4 we restrict our attention to Euclidean space, and show that the Nelson-Aalen estimator satisfies a functional central limit theorem. These results are applied in §5, to show how the Nelson-Aalen estimator can be used to develop various tests involving the dependence structure of the underlying distribution. The validity of bootstrap methods is established in each case.

Most of these results can be found in additional detail in [2] (the doctoral thesis of the first author), where estimation of the survival function is also considered.

# 2 Framework and Definitions

## 2.1 The Set-up

The framework we will be using is very similar to that used in [6], [8] and [9], except for some details. Let $T$ be a compact complete separable metric space and $\lambda$ a finite measure on $\mathcal{B}$, the Borel sets of $T$. A simple example is $T = [0,1]^d$, and we will return to this space in the sequel. Preserving the notation in [9], if $\mathcal{D}$ is an arbitrary class of sets, then $\mathcal{D}(u)$ will denote the class of finite unions of sets from $\mathcal{D}$. For $D$ an arbitrary subset of $T$, let $T_D$ be a countable dense subset of $D$. We will use '$\subset$' to indicate strict inclusion; moreover, $\overline{D}$ denotes the closure and $D^\circ$ denotes the interior of $D$. All processes will be indexed by sets belonging to an indexing collection $\mathcal{A}$:

**Definition 2.1** *A nonempty class $\mathcal{A}$ of compact, connected subsets of $T$ is called an* indexing collection *if it satisfies the following:*

1. *$\emptyset, T \in \mathcal{A}$, and for $A \in \mathcal{A}$, $A^\circ \neq A$ if $A \neq \emptyset$ or $T$.*

2. *$\mathcal{A}$ is closed under arbitrary intersections and if $A, B \in \mathcal{A}$ are nonempty, then $A \cap B$ is nonempty. Denote $\emptyset' = \cap_{A\in\mathcal{A}, A\neq\emptyset} A$. If $(A_i)$ is an increasing sequence in $\mathcal{A}$, then $\overline{\cup_i A_i} \in \mathcal{A}$.*

3. *The $\sigma$-algebra generated by $\mathcal{A}$, $\sigma(\mathcal{A}) = \mathcal{B}$.*

4. Separability from above:
*There exists an increasing sequence of countable subclasses $\mathcal{A}_n = \{A_1^n, A_2^n, ...\}$ of $\mathcal{A}$ closed under intersections and satisfying $\emptyset, \emptyset' \in \mathcal{A}_n(u)$. As well, there exists a sequence of functions $g_n : \mathcal{A} \to \mathcal{A}_n(u)$ such that*

   *(a) $g_n$ preserves arbitrary intersections and finite unions (i.e. $g_n(\cap_{A\in\mathcal{A}'} A) = \cap_{A\in\mathcal{A}'} g_n(A)$ for any $\mathcal{A}' \subseteq \mathcal{A}$, and if $\cup_{i=1}^k A_i = \cup_{j=1}^m A'_j$, then $\cup_{i=1}^k g_n(A_i) = \cup_{j=1}^m g_n(A'_j)$),*

   *(b) for each $A \in \mathcal{A}$, $A \subseteq (g_n(A))^\circ$,*

   *(c) $g_n(A) \subseteq g_m(A)$ if $n \geq m$,*

   *(d) for each $A \in \mathcal{A}$, $A = \cap_n g_n(A)$ ,*

   *(e) if $A, A' \in \mathcal{A}$ then for every $n$, $g_n(A) \cap A' \in \mathcal{A}$, and if $A' \in \mathcal{A}_n$ then $g_n(A) \cap A' \in \mathcal{A}_n$.*

   *(f) $g_n(\emptyset) = \emptyset$ $\forall n$.*

**Definition 2.2** *For $t \in T$, the 'past' of $t$ is the set $A_t := \cap_{A\in\mathcal{A},\ t\in A} A$ .*

We observe that the indexing collection $\mathcal{A}$ induces a partial order on $T$: $s \leq t$ if and only if $A_s \subseteq A_t$. In other words, $A_t = \{s : s \leq t\}$. We can now define the 'future' of $t \in T$ to be the set $E_t = \{s : s \geq t\}$. It will be assumed that $E_t$ is closed. In terms of this partial order, $\emptyset'$ is the minimal element of $T$.

The following definitions can be found in [9]: $\mathcal{C}$ is the class of all subsets of $T$ of the form

$$C = A \setminus B,\ A \in \mathcal{A},\ B \in \mathcal{A}(u).$$

$\mathcal{C}$ is closed under intersections and any set in $\mathcal{C}(u)$ may be expressed as a finite disjoint union of sets in $\mathcal{C}$; in other words, $\mathcal{C}$ is a semi-algebra.

For each $\mathcal{A}_k$, let $\mathcal{C}^\ell(\mathcal{A}_k)$ denote the *left-neighbourhoods* of $\mathcal{A}_k$: that is, all nonempty sets of the form $C = A \setminus \bigcup_{A' \in \mathcal{A}_k, A \not\subseteq A'} A', A \in \mathcal{A}_k$. It is clear that $\mathcal{C}^\ell(\mathcal{A}_k)$ partitions $T$, and as observed in [8], the sequence $(\mathcal{C}^\ell(\mathcal{A}_k))_k$ is a dissecting system for $T$.

**Example 2.3** As an example of the type of structure we will be using, consider $T = [0, 1]^2$. Then the class $\mathcal{A} = \{[0, t] : t \in T\}$ satisfies the conditions of Definition 2.1, and $\mathcal{C}(u)$ consists of all finite unions of disjoint rectangles of the form $(s, t]$ for $s, t \in T$.

For $T = [0, 1]^2$ and $A \in T$, $A$ will be called a *lower layer* if $[0, t] \subseteq A$, $\forall t \in A$. The class $\mathcal{A} = \{A \subseteq T : A \text{ is a lower layer}\}$ is also an indexing collection. In both cases we have that $A_t = [0, t]$ and $E_t = [t, 1]$, and both indexing collections induce the usual partial order in $[0, 1]^2$.

It is important to note that we can define other indexing collections that induce different partial orders. Example 2.4 in [8] illustrates this point.

We will require some extra assumptions. The first one is about the structure of $\mathcal{C}$, while the second one is concerned with the sub-semilattices $\mathcal{A}_k$:

**Assumption 2.4** *If $C \in \mathcal{C}$, then there exists a representation $C = A \setminus \cup_{i=1}^m D_i$ $(A, D_1, ..., D_m \in \mathcal{A})$ such that if $A' \in \mathcal{A}$, $A' \cap C = \emptyset$, then $A' \subseteq \cup_{i=1}^m D_i$. This is called a* maximal *representation of $C$.*

**Assumption 2.5** *For each $C \in \mathcal{C}^\ell(\mathcal{A}_k)$, there exists a point $t_{C^-} \in T$ such that $C \subseteq E_{t_{C^-}}$, and*

$$\lim_{k\to\infty} \sup_{C\in\mathcal{C}^\ell(\mathcal{A}_k)} \sup_{t\in C} d(t_{C^-}, t) = 0,$$

*where $d$ is the metric on $T$. Moreover, $t_{C^-}$ can be chosen so that $t_{C^-} \in \partial C$.*

The point $t_{C^-}$ acts as a lower bound on the points in $C$: $t_{C^-} \le t\ \forall t \in C$. For example, if $\mathcal{A} = \{[0, t] : t \in T\}$ and $C = (t, t']$, then $t_{C^-} = t$.

**Definition 2.6** *Let $(\Omega, \mathcal{F}, P)$ be any complete probability space. A* filtration *indexed by $\mathcal{A}$ is a class of complete sub-$\sigma$-fields of $\mathcal{F}$, $\{\mathcal{F}_A : A \in \mathcal{A}\}$ such that*

- *If $A \subseteq B$, then $\mathcal{F}_A \subseteq \mathcal{F}_B$ $\forall A, B \in \mathcal{A}$.*
- *$\mathcal{F}_{\cap A_i} = \bigcap \mathcal{F}_{A_i}$ for any decreasing sequence $(A_i)$ in $\mathcal{A}$ (Monotone outer-continuity).*

We continue to follow the definitions in [9]. If $B \in \mathcal{A}(u)$, then $\mathcal{F}^0_B = \bigvee_{A \in \mathcal{A},\ A \subseteq B} \mathcal{F}_A$. The $\sigma$-algebras $\{\mathcal{F}^0_B : B \in \mathcal{A}(u)\}$ are complete and increasing, but they may not be monotone outer-continuous. To avoid this problem, we can define $\mathcal{F}_B = \bigcap_n \mathcal{F}^0_{g_n(B)}$ for $B \in \mathcal{A}(u)$. For $C \in \mathcal{C}(u) \setminus \mathcal{A}$, let $\mathcal{G}^*_C = \bigvee_{B \in \mathcal{A}(u),\ B \cap C = \emptyset} \mathcal{F}_B$, and for $A \in \mathcal{A},\ A \neq \emptyset$, define $\mathcal{G}^*_A = \mathcal{F}_\emptyset$. Note that $\{\mathcal{G}^*_C\}$ is a decreasing family of $\sigma$-fields: if $C \subseteq C'$, then $\mathcal{G}^*_{C'} \subseteq \mathcal{G}^*_C$. Finally, if $C = A \setminus B$ is a maximal representation of $C$, then $\mathcal{G}^*_C = \mathcal{F}_B$. This was shown in [6].

**Definition 2.7** *(cf. [9]):*

- *A ($\mathcal{A}$ -indexed)* stochastic process *$X = \{X_A : A \in \mathcal{A}\}$ is a collection of random variables indexed by $\mathcal{A}$, and is said to be* adapted *if $X_A$ is $\mathcal{F}_A$-measurable for every $A \in \mathcal{A}$.*
- *$X$ is said to be* integrable *if $E[\mid X_A \mid] < \infty$ $\forall A \in \mathcal{A}$.*
- *A process $X : \mathcal{A} \to \mathbf{R}$ is* increasing *if for every $\omega \in \Omega$, the function $X.(\omega)$ can be extended to a finitely additive function on $\mathcal{C}$ satisfying $X_{\emptyset'}(\omega) = 0$ and $X_C(\omega) \geq 0$, $\forall C \in \mathcal{C}$, and such that $X.(\omega)$ is monotone outer-continuous (i.e. if $(A_n)$ is a decreasing sequence of sets in $\mathcal{A}(u)$ such that $\cap_n A_n \in \mathcal{A}(u)$, then $\lim_n X_{A_n}(\omega) = X_{\cap_n A_n}(\omega)$).* (It is clear that any trajectory of an increasing process can be extended to a measure on $\mathcal{B}$. Note that an increasing process is not necessarily adapted.)
- *An integrable process $M = \{M_A, A \in \mathcal{A}\}$ is called a* strong martingale *if it is adapted and for any $C \in \mathcal{C}$, $E[M_C|\mathcal{G}^*_C] = 0$. If the process $M$ is not adapted, it will be called a* pseudo-strong martingale.
- *A process $\overline{X}$ is called a $*$-* compensator *of the process $X$ if it is increasing and the difference $X - \overline{X}$ is a pseudo-strong martingale.*

Before moving on to adapted random sets, we should note that although we are restricting ourselves to a compact set $T$, all of the definitions and developments throughout this work can easily be expanded to a locally compact space using the structure found in Definition 2.1 of [8].

## 2.2 Clouds and $*$-clouds

We start by giving the definition of an adapted random set, a notion that was first introduced in [9]. Recall that a closed set $D$ is a *domain* if $D = \overline{D^\circ}$. As in [9], let $\mathcal{K}$ be the class of domains $D$ in $T$ such that $\lambda(\partial D) = 0$, where $\lambda$ is the measure defined at the beginning of Section 2.1, and let $\mathcal{L}$ be the class of open sets that are complements of sets in $\mathcal{K}$.

**Definition 2.8** *A random set $\eta : \Omega \to \mathcal{B}$ is an* adapted random set *if for any $t \in T$, $\{\omega : t \in \eta(\omega)\} \in \mathcal{F}_{A_t}$.*

- *An adapted random set $\rho$ taking its values in $\mathcal{L}$ is a* cloud*.*
- *An adapted random set $\xi$ taking its values in $\mathcal{K}$ is an* anti-cloud*.*

We now define a new kind of random set, but before doing so we need another assumption.

**Assumption 2.9** *Define $D_t = \overline{E_t^c}$. Then $D_t \in \mathcal{A}(u)\ \forall\ t \in T$. Also, if $(t_n)$ is a sequence such that $t \leq t_n\ \forall n$ and $t_n \downarrow t$, then $D_t = \bigcap D_{t_n}$.*

This is crucial to some later developments, since $E_t^\circ = T \setminus D_t$ and so we will have that $E_t^\circ \in \mathcal{C}$. We may interpret $E_t^\circ$ as the 'strict' future of $t$, and we shall define the 'wide' history of $t$ to be $\mathcal{F}_t^* := \mathcal{G}^*_{E_t^\circ} = \mathcal{F}_{D_t}$.

Now we are in position to define $*$-adapted sets:

**Definition 2.10** *A random set $\eta : \Omega \to \mathcal{B}$ is a $*$-*adapted random set *if for any $t \in T$, $\{\omega : t \in \eta(\omega)\} \in \mathcal{F}_t^*$.*

- *A $*$-adapted random set $\rho$ taking its values in $\mathcal{L}$ is a $*$-*cloud.
- *A $*$-adapted random set $\xi$ taking its values in $\mathcal{K}$ is a $*$-*anti-cloud.

The following theorems were proven in [9] for anti-clouds. We give the statements for $*$-anti-clouds; the proofs are analogous to their anti-cloud counterparts. For details, see [2]. We recall from the definition of an increasing process that if $X$ is increasing, it can be extended as a measure to $\mathcal{B}$. Therefore, $X_\xi(\omega) := X_{\xi(\omega)}(\omega)$ is well-defined for any anti-cloud or *-anti-cloud.

**Theorem 2.11** *Let $X$ be an increasing process (or the difference of two increasing processes) on $\mathcal{A}$ and suppose that $\xi$ is a $*$-anti-cloud. Then for any $A \in \mathcal{A}$, both $X_{A\cap\xi}$ and $X_{A\cap\partial\xi}$ are random variables.*

Now we deal with the filtered process $X^\xi_A := X_{\xi\cap A}$.

**Theorem 2.12** *Let $\xi$ be a $*$-anti-cloud and $X = Y - W$, where $Y$ and $W$ are increasing processes such that $Y_{\partial\xi} = W_{\partial\xi} = 0$ a.s.*

1. *If $X$ is a (pseudo)-strong martingale, then $X^\xi$ is a pseudo-strong martingale.*
2. *$X^\xi$ will not generally be adapted, even if $X$ is.*

# 3 Hazard Estimation and Observability

## 3.1 The Nelson-Aalen Estimator

The model is the same as the one in [8] and [9]. Assume that $(\Omega, \mathcal{F}, P)$ is a complete probability space. Let $Y : \Omega \to T$ be a $T$-valued random variable, and $\mu(B) = P\{Y \in B\}$ its distribution. The *survival function* associated with $Y$ is $S(t) = \mu(E_t)$. We assume that $\mu$ is absolutely continuous with respect to $\lambda$ ($\lambda$ as in §2) and denote by $\mu'$ the Radon-Nikodym derivative of $\mu$ with respect to $\lambda$ on the Borel sets of $T$. Finally, we will assume throughout this section that the indexing collection is $\mathcal{A} = \{A_t : t \in T\}$.

**Definition 3.1**
- *For $t \in T$, the hazard function of $Y$ is $h$ where*

$$h(t) = \frac{\mu'(t)}{S(t)}.$$

*If $S(t) = 0$, $h(t)$ is defined to be zero.*

- *The integrated hazard function of $Y$ is $H$ where*

$$H_A = \int_A h(u)\lambda(du) \text{ for any } A \in \mathcal{A}.$$

Let $N = \{N_A, A \in \mathcal{A}\} = \{I_{\{Y\in A\}}, A \in \mathcal{A}\}$ be the single jump process associated with $Y$ and $\mathcal{F}^Y = \{\mathcal{F}^Y_A, A \in \mathcal{A}(u)\}$ its minimal filtration: $\mathcal{F}^Y_A = \sigma\{N_B : B \in \mathcal{A}, B \subseteq A\} \cup \{\mathcal{P}_0\}$, where $\mathcal{P}_0$ is the class of $P$-null sets. This filtration is in fact outer-continuous, as was proved in [7]. $N$ is increasing, and it was proved in [8] that it has a $*$-compensator:

**Proposition 3.2** *([8], Proposition 2.9) The process $\overline{N}$ defined by*

$$\overline{N}_A = \int_{A\cap A_Y} \mu(E_u)^{-1}\mu(du) = \int_{A\cap A_Y} h(u)\lambda(du) = \int_A I_{\{Y\in E_u\}} h(u)\lambda(du)$$

*is a $*$-compensator of the process $N$ with respect to its minimal filtration, where $A_Y(\omega) = A_{Y(\omega)}$.*

Suppose we have a $T$-valued random variable $Y$ whose associated single jump process $N$ is adapted to a filtration $\mathcal{F}$ and an $\mathcal{F}$-$*$-cloud $\rho$ with corresponding $*$-anti-cloud $\xi = \rho^c$. The filtered jump process $N^\xi_\cdot := I_{\{Y\in\xi\cap\cdot\}}$ corresponds to observing occurrences of $Y$ only on the complement of the $*$-cloud; then we can say that $N$ has been filtered by the $*$-cloud $\rho$.

**Example 3.3** This example shows that our framework includes the usual bivariate censoring model, as presented in [11]. We assume that $T = [0.1]^2$ (or any bounded rectangle in $\mathbf{R}^2_+$) and let $\lambda$ denote Lebesgue measure. Let $Y = (Y_1, Y_2) \in [0,1]^2$ be a two-dimensional failure time on a probability space $(\Omega, \mathcal{F}, P)$, and suppose that $\mathcal{F}$ satisfies $\mathcal{F}^Y_t \subseteq \mathcal{F}_t$ for every $t \in [0,1]^2$. Let $t = (t_1, t_2) \in [0,1]^2$, $\mathcal{F}^{(1)}_{t_1} = \bigvee_{t_2} \mathcal{F}_{(t_1,t_2)}$, $\mathcal{F}^{(2)}_{t_2} = \bigvee_{t_1} \mathcal{F}_{(t_1,t_2)}$, and let $\tau = (\tau_1, \tau_2)$ be a two-dimensional censoring time, where $\tau_i$ is an $\mathcal{F}^{(i)}$-stopping time for $i = 1, 2$.

Suppose we can observe $Y \wedge \tau = (Y_1 \wedge \tau_1, Y_2 \wedge \tau_2)$ and $I_{\{Y_i \le \tau_i\}}$, $i = 1, 2$. Now we can express everything in terms of sets, letting $A_t = [0, t]$, and $\xi = [0, \tau]$. $\xi$ is a $*$-anti-cloud, since

$$\{t \in \xi\} = \{(t_1, t_2) \in [0, \tau_1] \times [0, \tau_2]\} = \{t_1 \le \tau_1\} \cap \{t_2 \le \tau_2\} \in \mathcal{F}^{(1)}_{t_1} \otimes \mathcal{F}^{(2)}_{t_2} \subseteq \mathcal{F}^*_t.$$

Finally, the counting process of censored times is

$$N^\xi(t) = N^\xi_{A_t} = I_{\{Y\wedge\tau\le t, Y_1\le\tau_1, Y_2\le\tau_2\}} = I_{\{Y\in A_t\cap\xi\}}.$$

□

In what follows, if $\xi : \Omega \to \mathcal{K}$, let $\mathcal{F}^\xi$ denote the minimal ($\mathcal{A}$-indexed) filtration with respect to which $\xi$ is a $*$-anti-cloud.

**Definition 3.4** *Let $Y$ be a $T$-valued random variable and let $\mathcal{F}^Y$ be the minimal filtration generated by its associated jump process $N$. Let $\mathcal{F}$ be a filtration such that $\mathcal{F}^Y_A \subseteq \mathcal{F}_A$ $\forall A \in \mathcal{A}(u)$ and let $\xi$ be an $\mathcal{F}$-$*$-anti-cloud. $\xi$ is*

1. *weakly independent of $Y$ if the $*$-compensator of $N$ with respect to $\mathcal{F}$ is the same as the $*$-compensator with respect to $\mathcal{F}^Y$;*
2. *independent of $Y$ if $\mathcal{F}^Y$ is independent of $\mathcal{F}^\xi$ and $\mathcal{F}_A = \mathcal{F}^Y_A \vee \mathcal{F}^\xi_A$, $\forall A \in \mathcal{A}$.*

The following lemma is an immediate consequence of Theorem 2.12.

**Lemma 3.5** *Suppose that $\xi$ is a $*$-anti-cloud, that $Y$ and $\xi$ are weakly independent and that the filtration $\mathcal{F}$ satisfies $\mathcal{F}_A = \mathcal{F}^Y_A \vee \mathcal{F}^\xi_A$, $\forall A \in \mathcal{A}$. If $N$ is the single-jump process associated with $Y$, then $(N - \overline{N})^\xi$ is a pseudo-strong ($\mathcal{F}$-)martingale.*

Now we are able to define the Nelson-Aalen estimator of the integrated hazard function using filtered data. We will use the same definitions as in [9].

Henceforth, we assume that we have a sequence of i.i.d. $T$-valued random variables $(Y_i)$ with the same distribution as $Y$, as well as a sequence $(\xi_i)$ of $*$-anti-clouds independent of the $Y_i's$. Define the following processes:

$$\begin{aligned} N^{(n)}_A &= \sum_{i=1}^n I_{\{Y_i \in A\}}, \\ Z_n(t) &= \sum_{i=1}^n I_{\{Y_i \in E_t\}} I_{\{t \in \xi_i\}}, \qquad (1) \\ \overline{N}^{(n)}_A &= \int_A Z_n(t) h(t) \lambda(dt), \\ {N^{(n)}}^\xi_A &= \sum_{i=1}^n I_{\{Y_i \in A \cap \xi_i\}}. \end{aligned}$$

By independence and Lemma 3.5, $\overline{N}^{(n)}$ is a $*$-compensator for ${N^{(n)}}^\xi$, and the process

$$M^{(n)}_\cdot = {N^{(n)}}^\xi_\cdot - \int_\cdot Z_n(t) h(t) \lambda(dt) \qquad (2)$$

is a pseudo-strong martingale with respect to $\mathcal{F}$, the minimal filtration generated by the sequences $(Y_i)$ and $(\xi_i)$. Since

$$N^{(n)^\xi}(dt) = Z_n(t)h(t)\lambda(dt) + M^{(n)}(dt),$$

regarding $M^{(n)}$ as noise, we come to a set-indexed version of the Nelson-Aalen estimator for $H_A$:

$$\hat{H}_A^{(n)} = \int_A \frac{N^{(n)^\xi}(dt)}{Z_n(t)} = \sum_{\{i:Y_i \in A \cap \xi_i\}} (Z_n(Y_i))^{-1}. \tag{3}$$

We observe that $\hat{H}_A^{(n)} - H_A = \int_A \frac{M^{(n)}(dt)}{Z_n(t)}$; the following is analogous to Proposition 4.5 of [9]; for details, see [2].

**Proposition 3.6** *$\hat{H}^{(n)} - H$ is a pseudo-strong martingale.*

## 3.2 The Observability Problem

The Nelson-Aalen estimator defined in (3) is identical to that introduced in [9], with the exception that in [9] it was assumed that $\xi$ is an anti-cloud. We will now explain the reason for incorporating the more general censoring mechanism ($*$-clouds) into the survival model.

**Example 3.7** Ideally, we would like to have all the information regarding the events in $\mathcal{F}_t^*$ at time $t$. Unfortunately, this is generally not possible in practice, and we have to settle for somewhat more limited information. Specifically, if $\xi$ is a fully observable anti-cloud, typically we may only be able to observe

$$\mathcal{H}_t = \sigma\{I_{\{Y \in A \cap \xi\}}, I_{\{s \in \xi\}} : \ A \subseteq D_t, \ s \in D_t\} \subseteq \mathcal{F}_t^*.$$

We need the event $I_{\{Y \in E_t\}} I_{\{t \in \xi\}}$ to be observable in order to construct the estimator for the integrated hazard function. But the problem is that $I_{\{Y \in E_t\}} I_{\{t \in \xi\}}$ is not $\mathcal{H}_t$-measurable, and so the estimator cannot be used. Indeed, suppose $t \in \xi$ and $D_t \cap \xi^c \neq \emptyset$. If $Y$ is not observed in $D_t \cap \xi$, it is impossible to know whether $Y \in D_t \cap \xi^c$ or $Y \in E_t$. This situation is illustrated in Figure 1 for $T \subset \mathbf{R}_+^2$.

To correct this, we can define $\xi^*(\omega) = \{t : D_t \cap \xi^c(\omega) = \emptyset\}$. $\xi^*$ is a $*$-anti-cloud, since $\{t \in \xi^*\} = \{D_t \cap \xi^c = \emptyset\} = \{D_t \subseteq \xi\} = \bigcap_{s \in D_t} \{s \in$ $\xi\} = \bigcap_{s \in T_{D_t}} \{s \in \xi\} \in \mathcal{H}_t \subseteq \mathcal{F}_t^*$, recalling that $T_{D_t}$ stands for a countable dense subset of $D_t$. Figure 2 pictures the $*$-anti-cloud $\xi^*$ corresponding to the anti-cloud $\xi$ in Figure 1.

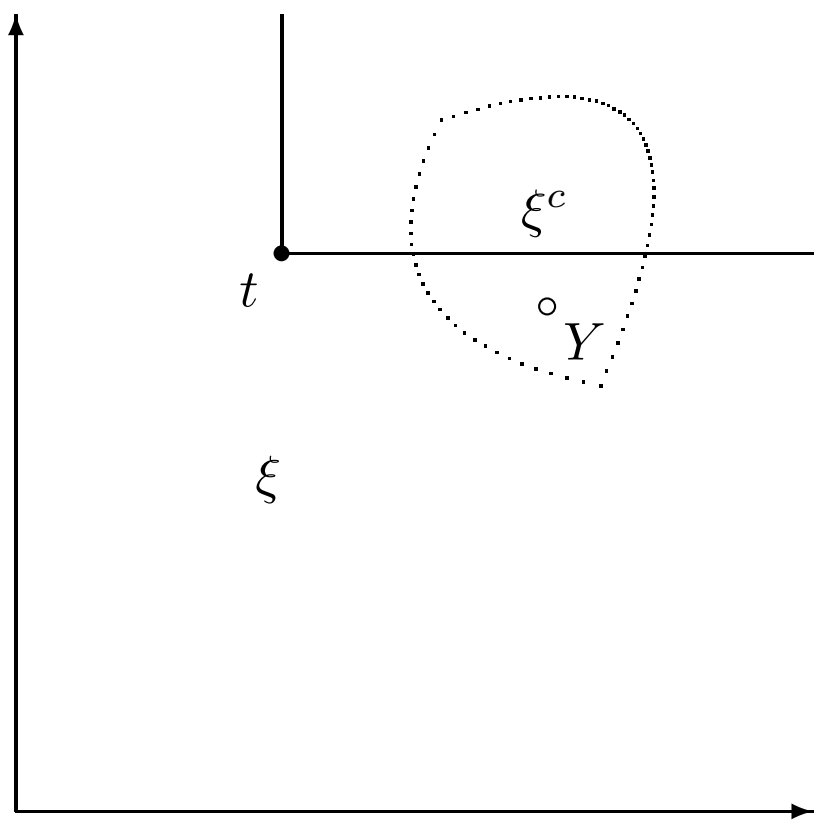

Figure 1: The observability problem: we don't know whether $Y \in E_t$.

Now consider $I_{\{Y \in E_t\}} I_{\{t \in \xi^*\}}$. We have that if $t \in \xi^*$, then $D_t \subseteq \xi$, and so

$$\begin{aligned} \{Y \in E_t\} \cap \{t \in \xi^*\} &= \{Y \in D_t\}^c \cap \{t \in \xi^*\} \\ &= \{Y \in D_t \cap \xi\}^c \cap \{t \in \xi^*\} \in \mathcal{H}_t. \end{aligned}$$

This means that $I_{\{Y \in E_t\}} I_{\{t \in \xi^*\}}$ is $\mathcal{H}_t$-measurable and hence observable, and now we are able to calculate the estimator.

Of course, there are instances where it is not possible to observe $I_{\{t \in \xi\}}$ whenever the observation $Y$ happens before time $t$; in other words, when the anti-cloud is not fully observable and $Y \in A_t$. In this case, the information available up to time $t$ consists of

$$\mathcal{H}_t = \sigma\{I_{\{Y \in A \cap \xi\}}, I_{\{s \in \xi \cap E_Y^c\}} : \ A \subseteq D_t, \ s \in D_t\} \subseteq \mathcal{F}_t^*;$$

nevertheless, we can still observe the event $I_{\{Y \in E_t\}} I_{\{t \in \xi^*\}}$: indeed, using the fact that $\{Y \in E_t\} = \bigcap_{s \in T_{D_t^\circ}} \{s \in E_Y^c\}$ and $\{t \in \xi^*\} = \bigcap_{s \in T_{D_t^\circ}} \{s \in \xi\}$, we have that

$$\{Y \in E_t\} \cap \{t \in \xi^*\} = \bigcap_{s \in T_{D_t^\circ}} (\{s \in E_Y^c\} \cap \{s \in \xi\}) = \bigcap_{s \in T_{D_t^\circ}} \{s \in E_Y^c \cap \xi\} \in \mathcal{H}_t.$$

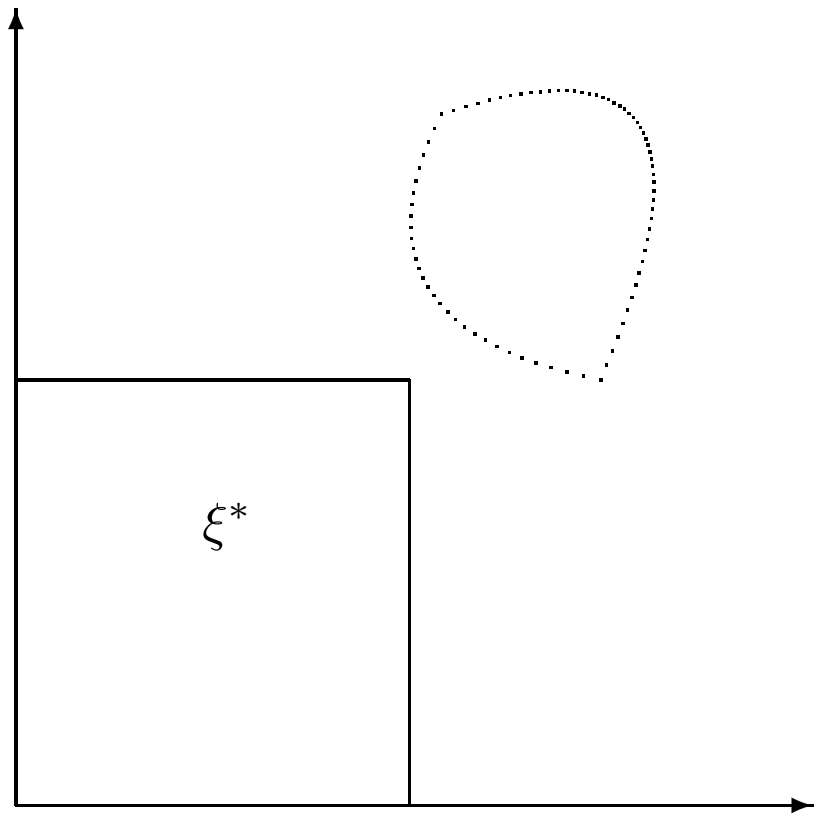

Figure 2: The resulting $*$-anti-cloud $\xi^*$.

This example demonstrates that an estimator can always be constructed, but in some sense is a worst-case scenario: we observe that the move from $\xi$ to $\xi^*$ as defined above may entail articially censoring observed values of $Y$ that lie in $\xi \setminus \xi^*$. Therefore, in practice, we should examine each information structure closely to determine whether the Nelson-Aalen estimator can be constructed using the cloud $\xi$ (cf. Example 3.10), and if not, how best to define an appropriate anti-cloud that censors as few observations as possible. This will be illustrated in Example 3.9. □

**Example 3.8** We now consider a generalization of Example 3.3. Again, assume that $T = [0,1]^2$ and $\lambda =$ Lebesgue measure. Let $Y = (Y_1, Y_2) \in [0,1]^2$ be a two-dimensional failure time on a probability space $(\Omega, \mathcal{F}, P)$, $t = (t_1, t_2) \in [0,1]^2$ with $\mathcal{F}_t^Y \subseteq \mathcal{F}_t$ for every $t \in [0,1]^2$ and $\mathcal{F}_{t_i}^{(i)}$ as in Example 3.3 for $i = 1, 2$; but instead of having a single two-dimensional censoring time, we will consider a finite sequence of $\mathcal{F}^{(i)}$-stopping times $0 = \eta_{i1} \leq \nu_{i1} \leq \eta_{i2} \leq \nu_{i2} \leq \ldots \leq \eta_{in_i} \leq \nu_{in_i}$ for $i = 1, 2$. Let $s = (s_1, s_2)$ and $u = (u_1, u_2)$ and suppose we can observe

$$\mathcal{H}_t = \{I_{\{Y_i \leq s_i, Y_i \in \bigcup_{j=1}^{n_i} [\eta_{ij}, \nu_{ij}]\}}, I_{\{Y_i \in (\nu_{ij}, \eta_{i,j+1}), \eta_{i,j+1} \leq u_i\}} : s, u \in A_t;\ i = 1, 2\}.$$

An illustration of this kind of situation can be seen in Figure 3. This scenario could arise in a laboratory, where test animals are under continuous observation during the day, but not at night.

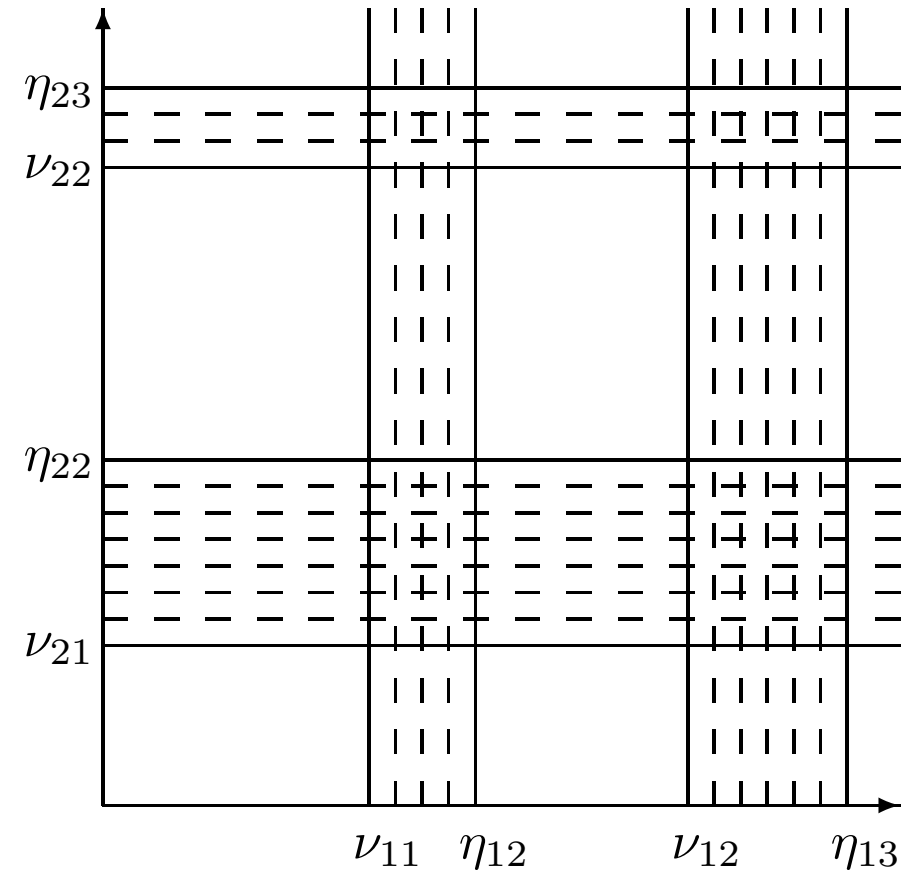

Figure 3: We can observe the failure times in the unshaded area.

Let $\xi^* = \{\bigcup_{j=1}^{n_1}[\eta_{1j}, \nu_{1j}] \times \bigcup_{j=1}^{n_2}[\eta_{2j}, \nu_{2j}]\}$. We have that

$$\begin{aligned}
\{t \in \xi^*\} &= \{(t_1, t_2) \in \bigcup_{j=1}^{n_1}[\eta_{1j}, \nu_{1j}] \times \bigcup_{j=1}^{n_2}[\eta_{2j}, \nu_{2j}]\} \\
&= \{t_1 \in \bigcup_{j=1}^{n_1}[\eta_{1j}, \nu_{1j}]\} \cap \{t_2 \in \bigcup_{j=1}^{n_2}[\eta_{2j}, \nu_{2j}]\} \\
&= (\bigcup_{j=1}^{n_1}\{t_1 \in [\eta_{1j}, \nu_{1j}]\}) \cap (\bigcup_{j=1}^{n_2}\{t_2 \in [\eta_{2j}, \nu_{2j}]\}) \\
&\in \mathcal{F}_{t_1}^{(1)} \otimes \mathcal{F}_{t_2}^{(2)} \subseteq \mathcal{F}_t^*,
\end{aligned}$$

so $\xi^*$ is a $*$-anti-cloud. Then $I_{\{Y \in E_t\}} I_{\{t \in \xi^*\}}$ can be observed and we get our estimator for the integrated hazard function. □

**Example 3.9** Once more, assume that $T = [0,1]^2$ and $\lambda$ = Lebesgue measure. Suppose now that we have a sequence of i.i.d. $T$-valued random variables $(Y_i) = (Y_{1,i}, Y_{2,i})$ with the same distribution as $Y$, as well as a sequence of i.i.d. $[0,1]$-valued random variables $\kappa_{1,i} < \kappa_{2,i}$ for $i = 1, 2, \ldots$. Define the filtration $\mathcal{F}_{t_1}^{\kappa_1,\kappa_2} = \sigma\{I_{\{\kappa_{1,i} \le s_1\}}, I_{\{\kappa_{2,i} \le s_1\}} : s_1 \le t_1\}$, and let $\mathcal{F}_{(t_1,t_2)} := \mathcal{F}_{(t_1,t_2)}^Y \vee \mathcal{F}_{t_1}^{\kappa_1,\kappa_2}$, as well as $\mathcal{F}_{t_1}^{(1)} = \bigvee_{t_2} \mathcal{F}_{(t_1,t_2)}$, $\mathcal{F}_{t_2}^{(2)} = \bigvee_{t_1} \mathcal{F}_{(t_1,t_2)}$. Let $\rho_i := \{(t_1, t_2) : \kappa_{1,i} < t_1 < \kappa_{2,i},\ t_1 < t_2 < t_1 + c\}$ represent the censored region, where $c$ denotes a constant. Then it is easy to see that $\rho_i$ is a sequence of $*$-clouds, since the $\kappa_{j,i}$ are stopping times. This kind of data structure can arise when, for example, $Y_1$ is the age of start of pregnancy and $Y_2$ is the age of onset of a disease such as tuberculosis. At $\kappa_1$, it is found that the test used for diagnosing the disease is dangerous for a pregnant woman, so once the pregnancy starts it becomes impossible to diagnose that disease for the next nine months. At time $\kappa_2$ though, a new safe test becomes available and the times of onset of the disease are no longer censored. See Figure 4.

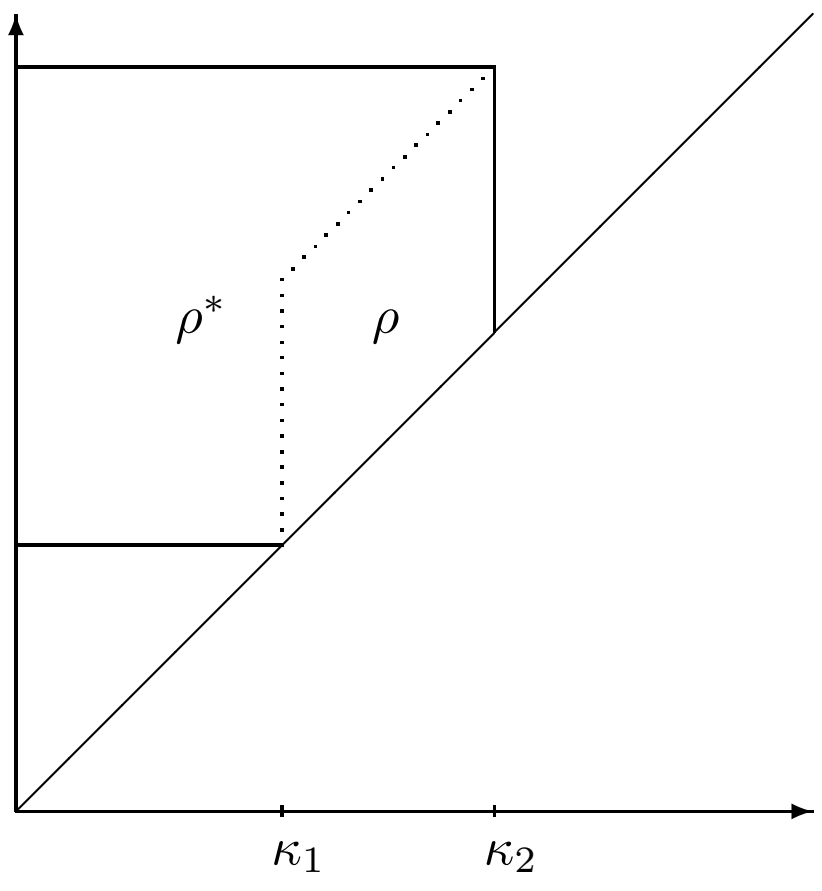

Figure 4: The regions $\rho$ and $\rho^*$.

In this case, we would censor the pairs $(Y_{1,i}, Y_{2,i})$ on the extended region $\xi_i^{*c} = \rho_i^* := \{(t_1, t_2) : \kappa_{1,i} < t_1 < \kappa_{2,i},\ t_1 < t_2 < \kappa_{2,i} + c\} \bigcup \{(t_1, t_2) : t_1 < \kappa_{1,i},\ \kappa_{1,i} < t_2 < \kappa_{2,i} + c\}$ in order to be able to observe the indicators $I_{\{Y \in E_t\}} I_{\{t \in \xi_i^*\}}$ needed to obtain the Nelson-Aalen estimator, as illustrated in Figure 4. Note that $\rho_i^*$ is also a sequence of $*$-clouds, since all the endpoints of the intervals in the definition of each $\rho_i^*$ are stopping times with respect to both $\mathcal{F}_{t_1}^{(1)}$ and $\mathcal{F}_{t_2}^{(2)}$; hence

$$\{(t_1, t_2) \in \rho_i^*\} \in \mathcal{F}_{t_1}^{(1)} \otimes \mathcal{F}_{t_2}^{(2)} \subseteq \mathcal{F}_{t_1}^{(1)} \vee \mathcal{F}_{t_2}^{(2)} = \mathcal{F}_{(t_1,t_2)}^*.$$

The censored region defined by $\rho^*$ is clearly much smaller than the $*$-clouds of Example 3.7, and would result in less lost data. □

**Example 3.10** There are examples that fit the original censoring model of [9] when $T = [0, 1]^2$ and $\lambda =$ Lebesgue measure. Suppose that $\xi$ is a cloud which is a lower layer and that we can observe

$$\mathcal{H}_t = \{I_{\{Y_i \le s_i\}}, I_{\{u \in \xi\}} : s, u \in A_t;\ i = 1, 2\}.$$

Since $\xi$ is a lower layer, $I_{\{Y \in E_t\}} I_{\{t \in \xi\}} = (1 - I_{\{Y_1 \le t_1\}})(1 - I_{\{Y_2 \le t_2\}}) I_{\{t \in \xi\}} \in \mathcal{H}_t$, and the cloud $\xi$ may be used in (3). An application of this model to a medical data set in given in [3]. □

# 4 A Functional Central Limit Theorem

In this section, for clarity of exposition we will assume that $T = [0,1]^2$ and $\lambda =$ Lebesgue measure, but all results can be extended to bounded rectangles in $\mathbf{R}^d$. Our objective is to prove a functional central limit theorem for the $T$-indexed process

$$U_t^{(n)} := \sqrt{n}(\hat{H}_{A_t}^{(n)} - H_{A_t}).$$

We will use the functional delta method, for which [13] is an excellent reference.

Throughout this section, we will make extensive use of the following notation. If $X$ is an arbitrary set, then the Banach space $l^\infty(X)$ is the set of all functions $f : X \to \mathbf{R}$ that are bounded uniformly and equipped with the norm $\|f\| = \sup_x |f(x)|$. $C[0,1]^2$ denotes the continuous functions, and $D[0,1]^2$ is the Banach space of all functions $f : [0,1]^2 \to \mathbf{R}$ continuous from the upper right quadrant and with limits from the other quadrants. Both $C[0,1]^2$ and $D[0,1]^2$ are equipped with the uniform norm. Finally, $BV_M[0,1]^2$ denotes the space of all functions in $D[0,1]^2$ with total variation bounded by $M$. Products of these spaces will always be equipped with a product norm.

The next lemma will be used in the proof of the functional central limit theorem. It is a two-dimensional version of Lemma 3.9.17 in [13], and since the proof is similar, it will be omitted.

**Lemma 4.1** *For each fixed $M$, the maps $\phi : l^\infty[0,1]^2 \times BV_M[0,1]^2 \to \mathbf{R}$ and $\psi : l^\infty[0,1]^2 \times BV_M[0,1]^2 \to D[0,1]^2$,*

$$\phi(A,B) = \int_{(0,1]} A dB, \qquad \psi(A,B)(t) = \int_{(0,t]} A dB \tag{4}$$

*are Hadamard-differentiable tangentially to $C[0,1]^2 \times D[0,1]^2$ at each $(A,B)$ in $l^\infty[0,1]^2 \times BV_M[0,1]^2$ such that $\int |dA| < \infty$, and the derivatives are given by*

$$\phi'_{A,B}(\alpha,\beta) = \int A d\beta + \int \alpha dB, \qquad \psi'_{A,B}(\alpha,\beta)(t) = \int_{(0,t]} A d\beta + \int_{(0,t]} \alpha dB,$$

*where $\int A d\beta$ is defined via the two-dimensional integration by parts formula found in Theorem 8.8 of [4] if $\beta$ is not of bounded variation.*

Before moving on to the central limit theorem, we have to make some assumptions that will allow us to apply the delta-method to our processes.

**Assumption 4.2** *$P(t \in \xi)$ is continuous in $t \in T$ and there exists $\epsilon > 0$ such that $P(t \in \xi) > \epsilon$ for every $t \in T$. For all $s, t \in T$, $P(s \in \xi) - P(s, t \in \xi) \leq K |s - t|$, where $K$ is a constant and $|\cdot|$ denotes the Euclidean norm.*

**Assumption 4.3** *The survival function $S$ satisfies a Lipschitz condition of order 1.*

**Assumption 4.4** *If $S(\tau, \tau) > 0$, then $g(\cdot) := (S(\cdot)P(\cdot \in \xi))^{-1}$ is of bounded variation on $[0, \tau]^2$.*

The following functional central limit theorem is the main result of this article.

**Theorem 4.5** *Under Assumptions 4.2, 4.4 and 4.3,*

$$\left( \frac{1}{\sqrt{n}} (Z_n(\cdot) - E(Z_n(\cdot))), \frac{1}{\sqrt{n}} ({N^{(n)}}^{\xi}_{A\cdot} - E({N^{(n)}}^{\xi}_{A\cdot})) \right) \Rightarrow (\Delta, \Gamma)$$

*where $\Delta$ and $\Gamma$ are tight Gaussian processes on $C[0,1]^2$ and $D[0,1]^2$ respectively. Moreover, $U^{(n)} \Rightarrow G$ in $l^\infty[0,\tau]^2$ for every $\tau$ such that $S(\tau, \tau) > 0$, where $G$ is a mean-zero Gaussian process such that for $t \in T$,*

$$G_t = \int_{[0,t]} \frac{d\Gamma(u)}{S(u)P(u \in \xi)} - \int_{[0,t]} \frac{\Delta(u)h(u)}{S(u)P(u \in \xi)} \, du.$$

*The first term in $G_t$ is defined by integration by parts.*

Before proceeding with the proof, we make a few observations and prove a lemma that will be required. The statement of the CLT is very similar to Example 3.9.19 in [13]. One of the differences, obviously, is that we are working on two dimensions instead of one. However, the major difference resides in our censoring mechanism and the information we are provided with; a consequence of this is that $Z_n$, the survivor function process, will not be in $D[0,1]^2$. This lack of sample path regularity necessitates Assumptions 4.2 and 4.3 above, which yield the following lemma and its corollary:

**Lemma 4.6** *Under Assumptions 4.2 and 4.3, the sequence of processes* $D_n(t) := \frac{1}{\sqrt{n}}[Z_n(t) - E(Z_n(t))]$ *converges in distribution in* $l^\infty[0,1]^2$ *to a tight Gaussian process* $\Delta$ *taking its values on* $C[0,1]^2$.

**Proof:** Let $W_t = I_{\{Y \in E_t\}}$ and $V_t = I_{\{t \in \xi\}}$. Note that Assumption 4.2 implies that $E\,|V_t - V_s| \leq 2K\,|t-s|$, since

$$E\,|V_t - V_s| = E|I_{\{t\in\xi\}} - I_{\{s\in\xi\}}| = P(t \in \xi, s \in \xi^c) + P(s \in \xi, t \in \xi^c).$$

Then we have that

$$\begin{aligned}
& E|I_{\{Y\in E_t\}}I_{\{t\in\xi\}} - I_{\{Y\in E_s\}}I_{\{s\in\xi\}}|) \\
&= E\,[E(|W_tV_t - W_tV_s|)|W_t, W_s] \\
&= E\,[|V_t - V_s|\,|W_t = W_s = 1]\,P(W_t = W_s = 1) \\
&+ E\,[V_t|W_t = 1, W_s = 0]\,P(W_t = 1, W_s = 0) \\
&+ E\,[V_s|W_t = 0, W_s = 1]\,P(W_t = 0, W_s = 1) \\
&\leq 2K\,|t-s|\,S(s \vee t) + P(t \in \xi)(S(t) - S(s \vee t)) \\
&+ P(s \in \xi)(S(s) - S(s \vee t)) \\
&\leq K^*\,|t-s|\,,
\end{aligned}$$

for some constant $K^*$, where the first inequality follows from the independence of the processes $V$ and $W$ and the observation at the beginning of this proof; the second is a consequence of Assumption 4.3. But by a nearly identical argument to that in Example 2.11.14 in [13] -we use two-dimensional blocks instead of closed intervals- this implies that the sequence of processes defined by

$$\frac{1}{\sqrt{n}}[\sum_{i=1}^{n} I_{\{Y_i\in E_t\}}I_{\{t\in\xi_i\}} - nS(t)P(t\in\xi)] = \frac{1}{\sqrt{n}}[Z_n(t) - E(Z_n(t))] = D_n(t)$$

converges in distribution to a tight Gaussian process $\Delta$ on $l^\infty[0,1]^2$. Continuity of $\Delta$ is a consequence of the fact that $E|D_n(s) - D_n(t)|^2 \leq K_1|s-t|$ for some constant $K_1 < \infty$ (cf. [13], pg. 41): from Assumptions 4.2 and 4.3, we have that

$$\begin{aligned}
E|D_n(s) - D_n(t)|^2 &= S(s)P(s\in\xi) + S(t)P(t\in\xi) - 2S(s\vee t)P(s,t\in\xi) \\
&\quad -[S(s)P(s\in\xi) - S(t)P(t\in\xi)]^2 \\
&= [S(s)P(s\in\xi) - S(s\vee t)P(s,t\in\xi)] + [S(t)P(t\in\xi)
\end{aligned}$$

$$
\begin{aligned}
& -S(s\vee t)P(s,t\in\xi)] - [S(s)P(s\in\xi) - S(t)P(t\in\xi)]^2 \\
= \;& [S(s) - S(s\vee t)]P(s\in\xi) + [S(t) - S(s\vee t)]P(t\in\xi) \\
& +S(s\vee t)[P(s\in\xi) - P(s,t\in\xi)] \\
& +S(s\vee t)[P(t\in\xi) - P(s,t\in\xi)] \\
& -[S(s)P(s\in\xi) - S(t)P(t\in\xi)]^2 \\
\leq \;& K_1\left|s-t\right|.
\end{aligned}
$$

□

**Corollary 4.7** *Under Assumptions 4.2 and 4.3*

$$
\sup_{t\in[0,1]^2}\left|\frac{Z_n(t)}{n} - S(t)P(t\in\xi)\right| \to_P 0
$$

*as* $n\to\infty$.

We are now ready to proceed with proof of our main result.

**Proof of Theorem 4.5:** Since the $\xi_i$'s are i.i.d., we have that the sequence of processes $(C_n)$ defined by

$$
C_n(t) := \frac{1}{\sqrt{n}}[\sum_{i=1}^{n} I_{\{Y_i\in A_t\cap\xi_i\}} - nP(Y\in A_t, Y\in\xi)]
$$

converges in distribution to a tight Gaussian process $\Gamma$ on $D[0,1]^2$. This follows simply from the CLT for empirical processes, since we are working with a subdistribution. Combined with Lemma 4.6, this means that

$$
\frac{1}{\sqrt{n}}[(Z_n(\cdot), N^{(n)}{}^{\xi}_{A.}) - (E(Z_n(\cdot)), E(N^{(n)}{}^{\xi}_{A.}))] \Rightarrow (\Delta,\Gamma)
$$

on $l^\infty[0,1]^2\times D[0,1]^2$, where $(\Delta,\Gamma)$ is a Gaussian process on $C[0,1]^2\times D[0,1]^2$. Recalling (3), we note that the estimator depends on the pair $(\frac{Z_n}{n}, \frac{N^{(n)\xi}}{n})$ through the maps

$$
(A,B) \longrightarrow (\frac{1}{A}, B) \longrightarrow \int \frac{1}{A}dB. \tag{5}
$$

It is a consequence of Lemma 4.1 that the map (5) is Hadamard-differentiable tangentially to $C[0,1]^2 \times D[0,1]^2$ on a domain of the type $\{(A,B) : \int |dB| \leq M, A \geq \epsilon\}$ for given $M$ and $\epsilon > 0$, at every point $(A,B)$ such that $1/A$ is of bounded variation. Whenever $t$ is restricted to an interval $[0,\tau]^2$ such that $S(\tau) > 0$, by Assumption 4.4 and Corollary 4.7, the pair $(\frac{Z_n}{n}, \frac{N^{(n)\xi}}{n})$ is contained in this domain with probability tending to 1 for $M \geq 1$ and $\epsilon$ sufficiently small. The derivative map is given by

$$(\alpha, \beta) \longmapsto \int (1/A) d\beta - \int (\alpha/A^2) dB.$$

Now we can apply the delta-method to conclude that

$$U_\cdot^{(n)} = \sqrt{n}(\hat{H}^{(n)}{}_{A_\cdot} - H_{A_\cdot}) \Rightarrow G_\cdot,$$

where

$$G_t = \int_{[0,t]} \frac{d\Gamma(u)}{S(u)P(u \in \xi)} - \int_{[0,t]} \frac{\Delta(u)}{(S(u)P(u \in \xi))^2} dP(Y \in A_u, Y \in \xi)$$

is again Gaussian. As in Lemma 4.1, the first term in the limiting process has to be defined by integration by parts, since $\Gamma$ may not be of bounded variation. Finally, we observe that $P(Y \in A_t, Y \in \xi) = \int_{A_t} P(u \in \xi) d\mu(u) = \int_{A_t} S(u)P(u \in \xi)h(u)du$.

The covariance structure for $G$ will be given later in this section. □

En route to finding the covariance structure for $G$, we need the following proposition.

**Proposition 4.8** *Under Assumptions 4.2 and 4.3, if $A$ is a Borel subset of $[0,\tau]^2$ where $S(\tau,\tau) > 0$, then as $n \to \infty$,*

$$\left[\sqrt{n}(\hat{H}_A^{(n)} - H_A)\right] - \left[\int_A \frac{1}{S(t)P(t \in \xi)} \cdot \frac{M^{(n)}(dt)}{\sqrt{n}}\right] \to_P 0.$$

**Proof:** For any Borel set $A \subseteq [0,\tau]^2$,

$$\begin{aligned} \sqrt{n}(\hat{H}_A^{(n)} - H_A) &= \sqrt{n} \int_A \frac{1}{Z_n(t)} M^{(n)}(dt) \\ &= \int_A \left( \frac{n}{Z_n(t)} - \frac{1}{S(t)P(t \in \xi)} \right) \frac{M^{(n)}(dt)}{\sqrt{n}} \qquad (6) \\ &+ \int_A \frac{1}{S(t)P(t \in \xi)} \cdot \frac{M^{(n)}(dt)}{\sqrt{n}}. \qquad (7) \end{aligned}$$

Using Corollary 4.7, we can show that (6) converges in probability to 0 with exactly the same argument as in the proof of Theorem 5.1 in [8]. □

**Remark 4.9** Proposition 4.8 gives us the central limit theorem for the finite-dimensional distributions of $\sqrt{n}(\hat{H}_A^{(n)} - H_A)$ over a more general class of Borel sets, since $\int_A \frac{1}{S(t)P(t\in\xi)} \cdot \frac{M^{(n)}(dt)}{\sqrt{n}}$ is the normalized sum of $n$ i.i.d. processes, as noted in [8]. This comes up short of giving us a functional CLT since we need to prove tightness. It is for this reason that we believe the use of the Delta-Method is a more elegant approach to this particular problem.

The next proposition is identical to its counterpart in [9] (Corollary 4.7).

**Proposition 4.10** *Let $g$ be a continuous function on $T$. Then $\int_\cdot g(t)M^{(n)}(dt)$ is a pseudo-strong martingale and for $C, D \in \mathcal{C}$,*

$$\begin{aligned}
&\mathrm{Cov}\left(\int_C g(s)M^{(n)}(ds), \int_D g(t)M^{(n)}(dt)\right) \\
&= n\left[\int_{C\cap D} g^2(t)S(t)P(t\in\xi)h(t)dt\right. \\
&+ \left.\int\int_{l(C,D)} g(s)g(t)\mu(E_s\cap E_t)P(s,t\in\xi)h(s)h(t)dsdt\right],
\end{aligned} \tag{8}$$

*where*

$$\begin{aligned}
l(C,D) &= \{(c,d)\in C\times D : c\in A_d^c \cap E_d^c\} \\
&= \{(c,d)\in C\times D : d\in A_c^c \cap E_c^c\}.
\end{aligned}$$

Now that we have Propositions 4.8 and 4.10, we are in position to give the covariance structure for the limiting process in Theorem 4.5.

**Lemma 4.11** *The covariance structure for the process $G$ in Theorem 4.5 is given, for $C, D \in \mathcal{C}$, by*

$$\begin{aligned}
\mathrm{Cov}(G_C, G_D) &= \int_{C\cap D} (S(t)P(t\in\xi))^{-1}h(t)dt \\
&+ \int\int_{l(C,D)} \frac{\mu(E_s\cap E_t)P(s,t\in\xi)}{S(s)S(t)P(s\in\xi)P(t\in\xi)} h(s)h(t)dsdt.
\end{aligned} \tag{9}$$

**Proof:** From Proposition 4.8, we have that the only term that contributes to the covariance is (7), and so (9) follows by an application of Corollary 4.10. □

# 5 Applications

For all our applications, we will assume that $T = [0,1]^2$, and that $\mathcal{F}_t = \mathcal{F}_{(t_1,t_2)}$ is trivial if either $t_1 = 0$ or $t_2 = 0$. We recall that $\mathcal{F}^{(1)}_{t_1} = \vee_{t_2}\mathcal{F}_{(t_1,t_2)}$ and $\mathcal{F}^{(2)}_{t_2} = \vee_{t_1}\mathcal{F}_{(t_1,t_2)}$. Assumptions 4.2, 4.3 and 4.4 will be assumed to hold throughout this section.

## 5.1 Test of Independence

We now concern ourselves with the construction of a test of independence of the components $(Y_1, Y_2)$ of the random vector $Y \in T$. We will follow closely the development of the ideas presented in [11].

As noted in [8] and [11], when $Y_1$ and $Y_2$ are independent, the cumulative hazard is the product of the marginal hazards. Therefore, we have to be able to estimate the marginal hazards in order to obtain a test of independence. We will review the four examples given in §3.2 in order to illustrate how this is done.

• Example 3.7: Given $\xi$ and the corresponding $*$-anti-cloud $\xi^*$, define

$$\begin{aligned} \kappa_1 &:= \sup\{t_1 : (t_1, 0) \in \xi^*\} \\ \kappa_2 &:= \sup\{t_2 : (0, t_2) \in \xi^*\}, \end{aligned}$$

see Figure 5; note that $\kappa_1$ is an $\mathcal{F}^{(1)}$-stopping time: for any $t_1 \in [0,1]$,

$$\{\kappa_1 \le t_1\} = \{t_1 < \kappa_1\}^c = \{(t_1, 0) \in \xi_i^*\}^c \in \mathcal{F}^*_{(t_1,0)} = \mathcal{F}^{(1)}_{t_1}.$$

Similarly, we can show that $\kappa_2$ is an $\mathcal{F}^{(2)}$-stopping time. Taken separately, $Y_j$ is censored on the intervals $\xi_j^c := (\kappa_j, 1]$ for $j = 1,2$.

• Example 3.8: This example is illustrated in Figure 3. In this case, since $\nu_{ij}, \mu_{ij}$ are $\mathcal{F}^{(i)}$-stopping times, $i = 1,2;, j = 1,2,...$, $Y_i$ is censored on the $\mathcal{F}^{(i)}$-adapted clouds $\cup_j(\nu_{ij}, \mu_{ij}), i = 1,2$.

• Example 3.9: Referring to Figure 4, taken separately, $Y_1$ is not censored at all, since it is always possible to determine whether the individual is pregnant or not, but $Y_2$ is censored on the interval $(\kappa_{1,i}, \kappa_{2,i} + c)$.

• Example 3.10: In this structure, $Y_1$ is censored on the right by the $\mathcal{F}^{(1)}$-stopping time

$$\kappa_1 = \inf\{s : (s, 0) \notin \xi\}.$$

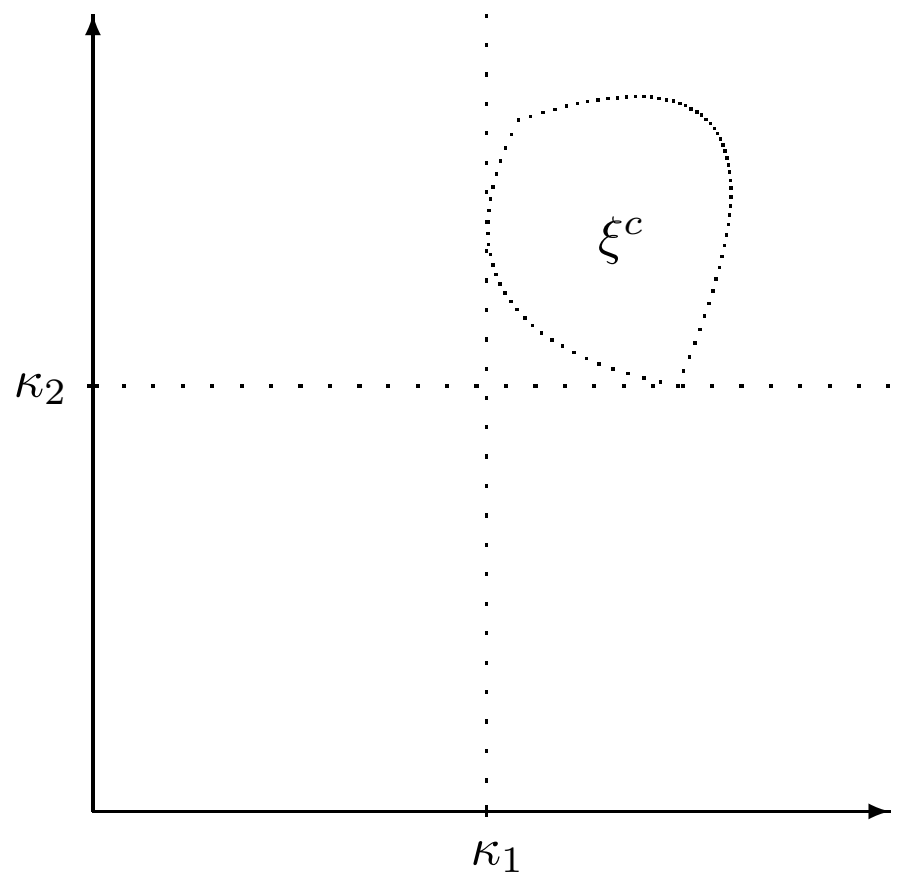

Figure 5: The stopping times $\kappa_1$ and $\kappa_2$

Likewise, $Y_2$ is censored on the right by the $\mathcal{F}^{(2)}$-stopping time

$$\kappa_2 = \inf\{u : (0, u) \notin \xi\}.$$

In each of these examples, we end up with the same type of one and two-dimensional structure: the pair $(Y_1, Y_2)$ is filtered on an $\mathcal{F}$-$*$-cloud, and each $Y_j$ is filtered by an $\mathcal{F}^{(j)}$-$*$-cloud $\xi_j^c$, a union of random open intervals with endpoints that are $\mathcal{F}^{(j)}$-stopping times for $j = 1, 2$. We would like to note that in classical one-dimensional problems, the intervals are usually taken to be left-open and right-closed in order to ensure predictability, but since we are assuming that $P(Y_{\partial\xi} = 0) = 1$, the endpoints of the intervals can actually be ignored.

Now we are in position to define the analogous one-dimensional processes

$$\begin{aligned}
N_1^{(n)^{\xi_1}}(t_1) &= \sum_{i=1}^{n} I_{\{Y_{1,i} \le t_1\}} I_{\{Y_{1,i} \in \xi_{1,i}\}} \\
N_2^{(n)^{\xi_2}}(t_2) &= \sum_{i=1}^{n} I_{\{Y_{2,i} \le t_2\}} I_{\{Y_{2,i} \in \xi_{2,i}\}} \\
Z_{n,1}(t_1) &= \sum_{i=1}^{n} I_{\{Y_{1,i} \ge t_1\}} I_{\{t_1 \in \xi_{1,i}\}} \\
Z_{n,2}(t_2) &= \sum_{i=1}^{n} I_{\{Y_{2,i} \ge t_2\}} I_{\{t_2 \in \xi_{2,i}\}},
\end{aligned}$$

as well as the processes $C_{n,j}(t_j) = n^{-1/2}[N_j^{(n)^{\xi_j}}(t_j) - E(N_j^{(n)^{\xi_j}}(t_j))]$ and $D_{n,j}(t_j) = n^{-1/2}[Z_{n,j}(t_j) - E(Z_{n,j}(t_j))]$ for $j = 1, 2$. We will also make use of the one-dimensional analogues of the process $M_\cdot^{(n)}$:

$$M_j^{(n)}(\cdot) = N_j^{(n)^{\xi_j}}(\cdot) - \int_\cdot Z_{n,j}(u_j)h(u_j)\lambda(du_j), \qquad j = 1, 2.$$

The cumulative marginal hazards are estimated by

$$\hat{H}_j^{(n)}(t_j) = \int_{[0,t_j]} \frac{N_j^{(n)^{\xi_j}}(du_j)}{Z_{n,j}(u_j)}, \qquad j = 1, 2.$$

**Lemma 5.1** $n^{1/2}(\hat{H}^{(n)} - H, \hat{H}_1^{(n)} - H_1, \hat{H}_2^{(n)} - H_2) \Rightarrow (G, G_1, G_2)$ *in the space* $l^\infty([0, \tau_1] \times [0, \tau_2]) \times l^\infty[0, \tau_1] \times l^\infty[0, \tau_2]$ *for every* $\tau = (\tau_1, \tau_2)$ *such that* $S(\tau) > 0$, *where the vector* $(G, G_1, G_2)$ *is jointly Gaussian with continuous sample paths.*

**Proof:** We start very similarly to the proof of Lemma 4.1 in [11]. Recall that $D_n(t) := n^{-1/2}[Z_n(t) - E(Z_n(t))]$ and $C_n(t) := n^{-1/2}[N^{(n)}{}^{\xi}_{A_t} - E(N^{(n)}{}^{\xi}_{A_t})]$. We already saw in the proof of Theorem 4.5 that $(D_n, C_n)$ converges to $(\Delta, \Gamma)$ on $l^\infty[0, 1]^2 \times D[0, 1]^2$. We also know that $(D_{n,j}, C_{n,j})$ converges to a joint Gaussian process $(\Delta_j, \Gamma_j), j = 1, 2$, on $l^\infty[0, 1] \times D[0, 1]$, since the same arguments for the convergence of $(D_n, C_n)$ still apply. Therefore, the joint sequence $(D_n, D_{n,1}, D_{n,2}, C_n, C_{n,1}, C_{n,2})$ is tight in the product of topologies of the space $l^\infty[0, 1]^2 \times (l^\infty[0, 1])^2 \times D[0, 1]^2 \times (D[0, 1])^2$. This, coupled with the fact that its finite-dimensional distributions converge to those of a Gaussian process with continuous sample paths, give the convergence of the sequence.

Next, note that both the estimators for the cumulative marginal hazards depend on the pairs $(\frac{Z_{n,j}}{n}, \frac{N_j^{(n)^{\xi_j}}}{n})$, $j = 1, 2$, through the map (5), where the integral is defined on $\mathbf{R}$ instead of $\mathbf{R}^2$. Using the same arguments as in Theorem 4.5,we can show that $n^{1/2}(\hat{H}_j^{(n)} - H_j) \Rightarrow G_j$ in $D[0, \tau_j]$ for every $\tau_j$ such that $S_j(\tau_j) > 0$, where $G_j$ are mean-zero Gaussian processes for $j = 1, 2$. As a consequence of this fact and the joint convergence of $(D_n, D_{n,1}, D_{n,2}, C_n, C_{n,1}, C_{n,2})$, we get the desired result. $\square$

As was the case in [11], in order to construct our test of independence between $(Y_{1,i})$ and $(Y_{2,i})$, $i = 1, 2, \ldots$, we will take the difference between the two-dimensional Nelson-Aalen estimator and the estimator under the hypothesis of independence. Define $V_n := n^{1/2}(\hat{H}^{(n)} - \hat{H}_1^{(n)}\hat{H}_2^{(n)})$ on $[0, \tau]$, where $\tau = (\tau_1, \tau_2)$. We could calculate the variance of $V_n$, but as Pons remarked in [12], its practical estimation is not easy and therefore not very useful. That is why we prefer to use a bootstrap test based on $V_n$, and our efforts will now turn towards that goal.

We shall use the delta-method for the bootstrap, and for that we need to verify that the classes of functions needed are Donsker classes. Let $X_1, \ldots, X_n$ be a sample of random elements in a measurable space $(\mathcal{X}, \mathcal{A})$ with distribution $P$, let $\mathcal{G}$ be a collection of measurable functions $g : \mathcal{X} \to \mathbf{R}$ and let $P_n$ be the empirical measure of the $X_i$, which induces a map from $\mathcal{G}$ to $\mathbf{R}$ defined by $P_n g := n^{-1} \sum_{i=1}^{n} g(X_i)$.

Ivanoff and Merzbach observed in [9] that an anti-cloud is measurable as a random closed set as a consequence of part 2 of their Corollary 3.10; it is easy to see that the same is true of $*$-anti-clouds. Now assume that we can fully observe the pair $(Y_i, \xi_i)$. Then we can take $(Y_i, \xi_i)$ to be the random elements $X_i$ in the previous definition. Let $P$ be the joint distribution of the pair $(Y_i, \xi_i)$, let $P_n$ denote its empirical measure, $\tilde{P}_n$ the bootstrap empirical distribution and $\tilde{G}_n$ the bootstrap empirical process defined by $\tilde{G}_n := \sqrt{n}(\tilde{P}_n - P_n)$. For a more complete description of these elements, as well as for the exact statement of the delta-method for bootstrap in probability -which we will use a few lines ahead- refer to Sections 3.6 and 3.9 of [13], respectively.

Recalling that $T = [0, 1]^2$, let $\mathcal{G}^1 = \{g_t^1 : t \in T\}$ and $\mathcal{G}^2 = \{g_t^2 : t \in T\}$ be defined, respectively, by $g_t^1(Y_i, \xi_i) := I_{\{Y_i \in E_t\}} I_{\{t \in \xi_i\}}$ and $g_t^2(Y_i, \xi_i) := I_{\{Y_i \in A_t\}} I_{\{Y_i \in \xi_i\}}$. Then both $\mathcal{G}^1$ and $\mathcal{G}^2$ are Donsker classes, since both were shown to be convergent in distribution to tight Gaussian limits on $l^\infty[0, 1]^2$ in the proof of Theorem 4.5. We can apply Theorem 2.10.6 in [13] to find that the pair $(\mathcal{G}^1, \mathcal{G}^2)$ is also a Donsker class, since each one is uniformly bounded. It is also clear that $\mathcal{G}^1$ and $\mathcal{G}^2$ have finite envelope functions because both classes are comprised exclusively of indicator functions. Hence, by Theorem 3.6.1 in [13], the conditions (3.9.9) of [13] are met. Furthermore, we have shown that the map (5) is Hadamard-differentiable tangentially to $C[0, 1]^2 \times D[0, 1]^2$ on a certain domain. Then we only need to apply Theorem 3.9.11 in [13], the delta-method for bootstrap in probability, to get the conditional convergence given $(Y_i, \xi_i)$ of the bootstrapped Nelson-Aalen estimator. Similar arguments lead to the same conclusion for the one-dimensional bootstrapped Nelson-Aalen estimators.

The above is summarized in the next lemma, which is an analogue of Lemma 5.1 for the bootstrapped estimators.

**Lemma 5.2** *Assume we can observe the pairs* $(Y_i, \xi_i)$ *for* $i = 1, 2, \ldots,$ *and let* $(\tilde{Y}_i, \tilde{\xi}_i)$ *denote the bootstrap sample. Then if* $\tilde{H}^{(n)}, \tilde{H}_1^{(n)}, \tilde{H}_2^{(n)}$ *represent the bootstrapped versions of the joint and marginal Nelson-Aalen estimators, we have that* $n^{1/2}(\tilde{H}^{(n)} - \hat{H}^{(n)}, \tilde{H}_1^{(n)} - \hat{H}_1^{(n)}, \tilde{H}_2^{(n)} - \hat{H}_2^{(n)})$ *converges conditionally given* $(Y_i, \xi_i)$ *to the Gaussian process* $(G, G_1, G_2)$ *of Lemma 5.1.*

Now that we have seen that the bootstrapping is working correctly under the assumption of complete observability of the pairs $(Y_i, \xi_i)$, we note that the bootstrapped version of the statistics only involve the product of indicator functions $q_i = I_{\{Y_i \in A_t\}} I_{\{Y_i \in \xi_i\}}$ and $r_i = I_{\{Y_i \in E_t\}} I_{\{t \in \xi_i\}}$. Therefore, we can drop the assumption of complete observability, since those indicators can still be constructed because the events involved are observable under the survival model. Therefore, we can in fact bootstrap directly from the observed values $q_i, r_i$ and come up with the same estimators.

We can deduce from Lemma 5.1 that on $[0, \tau_1] \times [0, \tau_2]$,

$$(n^{1/2}[\hat{H}^{(n)}(t_1, t_2) - \hat{H}_1^{(n)}(t_1)\hat{H}_2^{(n)}(t_2) - (H(t_1, t_2) - H_1(t_1)H_2(t_2))])$$

converges in distribution to a mean-zero Gaussian process; indeed, we can rewrite $V_n = n^{1/2}(\hat{H}^{(n)} - \hat{H}_1^{(n)}\hat{H}_2^{(n)})$ as

$$\begin{aligned} V_n &= n^{1/2}(\hat{H}^{(n)} - H) - n^{1/2}(\hat{H}_1^{(n)} - H_1)\hat{H}_2^{(n)} - n^{1/2}(\hat{H}_2^{(n)} - H_2)H_1 \\ &+ n^{1/2}(H - H_1 H_2), \end{aligned}$$

hence

$$\begin{aligned} n^{1/2}\left(\hat{H}^{(n)} - \hat{H}_1^{(n)}\hat{H}_2^{(n)} - (H - H_1H_2)\right) &= n^{1/2}(\hat{H}^{(n)} - H) \\ &- n^{1/2}(\hat{H}_1^{(n)} - H_1)(\hat{H}_2^{(n)} - H_2) \\ &- n^{1/2}(\hat{H}_1^{(n)} - H_1)H_2 \\ &- n^{1/2}(\hat{H}_2^{(n)} - H_2)H_1, \end{aligned}$$

which converges weakly to $G - H_2G_1 - H_1G_2$. Lemma 5.2 guarantees us that for the bootstrapped estimators,

$$\sup_{t \in [0,\tau]} \left| n^{1/2}[\tilde{H}^{(n)}(t_1, t_2) - \tilde{H}_1^{(n)}(t_1)\tilde{H}_2^{(n)}(t_2) - (\hat{H}^{(n)}(t_1, t_2) - \hat{H}_1^{(n)}(t_1)\hat{H}_2^{(n)}(t_2))] \right| \quad (10)$$

converges in distribution to the sup of the absolute value of that same Gaussian process. Then we can bootstrap from $(q_i, r_i)$ jointly and find $c_\alpha$ so that $(1-\alpha)100\%$ of the absolute values of (10) fall below $c_\alpha$, which would give us $(1-\alpha)100\%$ uniform confidence bands for the difference $H - H_1 H_2$.

Finally, for our test of independence $H_\circ : H = H_1 H_2$ vs. $H_\circ : H \neq H_1 H_2$, we reject $H_\circ$ if $\sup_{t \in [0,\tau]} n^{1/2} \left| \hat{H}^{(n)}(t_1, t_2) - \hat{H}_1^{(n)}(t_1) \hat{H}_2^{(n)}(t_2) \right| > c_\alpha$.

## 5.2 Test of Hazard Rate Order

Let $X_1, \ldots, X_n$ and $Y_1, \ldots, Y_m$ be independent random samples from bivariate distributions $F$ and $G$ respectively, and suppose that there is a common censoring mechanism in the form of a sequence of independent $*$-clouds. We are interested in testing whether the distributions are equal on a particular set or, more generally, on a suitable class of sets, against the alternative that there is a difference in the hazard rates.

We start with the case where we look at the hazards on a fixed set $A$. The null hypothesis is $H_\circ : F = G$, and we can test it either against a single-sided alternative such as $H_A^F < H_A^G$, or the double-sided alternative $H_A^F \neq H_A^G$. Tests using a one-sided alternative leads us to tests of hazard rate order such as $H_\circ : F = G$ vs. $H_1 : h^F < h^G$ on $A$, since this alternative would imply $H_A^F < H_A^G$. It is natural to consider the difference between the Nelson-Aalen estimators $\hat{H}^F$ and $\hat{H}^G$ of the integrated hazards, which leads to the test statistic

$$W_A^{(N)} = \sqrt{\frac{nm}{N}} \left( \hat{H}_A^{(n)F} - \hat{H}_A^{(m)G} \right),$$

where $N = n + m$.

By Remark 4.9, we still have that both $U_A^{(n)F} = \sqrt{n}(\hat{H}_A^{(n)F} - H_A^F)$ and $U_A^{(m)G} = \sqrt{m}(\hat{H}_A^{(m)G} - H_A^G)$ converge to independent mean-zero Gaussian limits $U_A^F$ and $U_A^G$ respectively. We can write

$$W_A^{(N)} = \sqrt{\frac{m}{N}} U_A^{(n)F} - \sqrt{\frac{n}{N}} U_A^{(m)G} + \sqrt{\frac{nm}{N}} \left( H_A^F - H_A^G \right).$$

Under the null hypothesis of equality of the cumulative hazards on the set $A$, and assuming $\frac{n}{n+m} \to \lambda$ as $n, m \to \infty$, we get that

$$W_A^{(N)} \Rightarrow \sqrt{1-\lambda} U_A^F - \sqrt{\lambda} U_A^G.$$

The limit variable has the same distribution as $U_A^F$ under the null hypothesis. Then we have a test of asymptotic level $\alpha$ for a two-sided alternative if we reject the null hypothesis whenever $\left|W_A^{(N)}\right| > w^{(N)}$, where we choose $w^{(N)}$ such that $w^{(N)} \to w^F = \inf\left\{t : P(\left|U_A^F\right| > t) \le \alpha\right\}$. A test with a one-sided alternative would be handled similarly; we just have to remove the absolute values.

We can also consider testing for the hazard rate over a whole class, as long as it is Donsker in order to preserve the Gaussian limits. An example would be to take the class of rectangles $\mathcal{A} = \{A_z = [0, z] : z \in [0, \tau]\}$, where $[0, \tau] = [0, \tau_1] \times [0, \tau_2]$ for $\tau = (\tau_1, \tau_2)$. All the previous remarks about $W^{(N)}$ are still valid, but now use $\sup_{z\in[0,\tau]} \left|W_{A_z}^{(N)}\right|$ and $\sup_{z\in[0,\tau]} W_{A_z}^{(N)}$ as the test statistics, as well as choosing $w^{(N)} \to w^F = \inf\left\{t : P(\sup_{z\in[0,\tau]} \left|U_{A_z}^F\right| > t) \le \alpha\right\}$ and $w^{(N)} \to w^F = \inf\left\{t : P(\sup_{z\in[0,\tau]} U_{A_z}^F > t) \le \alpha\right\}$ for two-sided and one-sided alternatives, respectively.

Assuming we can fully observe the pairs $(X_i, \xi_i)$, $(Y_j, \xi_j')$ for $i = 1, \dots, n$ and $j = 1, \dots, m$, we determine the appropriate critical value for $w^{(N)}$ by bootstrapping from the pooled sample $S^{(N)} = (X_1, \xi_1), \dots, (X_n, \xi_n), (Y_1, \xi_1'), \dots, (Y_m, \xi_m')$, where we choose any of these observations with probability $1/N$, as was done in Section 3.7.2 in [13]. Define $J = \lambda F + (1 - \lambda)G$, and note that under the null hypothesis, $H^J = H^F$. We assign the first $n$ elements from the resampling to $F$, and the rest to $G$. Let $\ddot{H}^{(n,N)F}$ and $\ddot{H}^{(m,N)G}$ be the Nelson-Aalen estimators for the pooled sample assigned to $F$ and $G$ respectively, and $\ddot{H}^{(N)}$ the estimator for the complete pooled sample. Set

$$\ddot{W}^{(N)} = \sqrt{\frac{nm}{N}} \left(\ddot{H}^{(n,N)F} - \ddot{H}^{(m,N)G}\right).$$

Next, note that by Hadamard-differentiability and Theorem 3.7.6 in [13], we know that $\sqrt{n}\left(\ddot{H}^{(n,N)F} - \ddot{H}^{(N)}\right) \Rightarrow U_1^J$ and $\sqrt{m}\left(\ddot{H}^{(m,N)G} - \ddot{H}^{(N)}\right) \Rightarrow U_2^J$, where $U_1^J$ and $U_2^J$ are independent mean-zero Gaussian processes with covariance as defined in Lemma 4.11, with $S$ equal to the survival function of $J$ and $h = J'/S$. Then if we let $\lambda_N = n/N$,

$$\ddot{W}^{(N)} = \sqrt{n(1 - \lambda_N)} \left(\ddot{H}^{(n,N)F} - \ddot{H}^{(N)}\right) - \sqrt{m\lambda_N} \left(\ddot{H}^{(m,N)G} - \ddot{H}^{(N)}\right),$$

which converges in distribution to $\sqrt{1 - \lambda}U_1^J - \sqrt{\lambda}U_2^J$, also a Gaussian process on $[0, \tau]$ equal in distribution to $U_1^J$. Then, as remarked in Section 3.7.2 of [13], we can use

$$\ddot{w}_N = \inf\left\{t : P(\sup_{z\in[0,\tau]} \ddot{W}_{A_z}^{(N)} > t) \le \alpha\right\}$$

as critical values for a one-sided test, and

$$\ddot{w}_N = \inf\left\{ t : P(\sup_{z\in[0,\tau]} |\ddot{W}^{(N)}_{A_z}| > t) \leq \alpha \right\}$$

for the two-sided test.

Once again, as we did following the statement of Lemma 5.2, we can see that the bootstrapped version of our test statistic only uses the functions $q_i = I_{\{X_i\in A_t\}} I_{\{X_i\in\xi_i\}}$, $r_i = I_{\{X_i\in E_t\}} I_{\{t\in\xi_i\}}$, $q'_j = I_{\{Y_j\in A_t\}} I_{\{Y_j\in\xi'_j\}}$ and $r'_j = I_{\{Y_j\in E_t\}} I_{\{t\in\xi'_j\}}$; then we can safely drop the assumption of complete observability, given that even without it we can still construct the appropriate estimators. Thus, the test can be performed by bootstrapping directly from the recorded values $q_i, r_i, q'_j, r'_j$ and following the procedure described above.

## 5.3 The Hazard of a Copula

Suppose we have a sequence $(X_i, Y_i)$ of i.i.d. bivariate random vectors with continuous distribution $J$ and continuous, strictly increasing marginals $F$ and $G$, and let $H$ denote the integrated hazard function. Our goal is to estimate the hazard of the copula $C$ associated with $J$.

As usual, each pair $(X_i, Y_i)$ will be censored by an $\mathcal{F}$-$*$-cloud $\xi_i^c$. Let $C(p,q) = J(F^{-1}(p), G^{-1}(q))$, the copula function for $J$. and let $\overline{C}$ be the survival function for the copula: $\overline{C}(p,q) = 1 - p - q + C(p,q)$. It is straightforward to verify that the integrated hazard of the copula, $H^C$, satisfies

$$H^C(p,q) = H(F^{-1}(p), G^{-1}(q)), 0 \leq p, q \leq 1. \tag{11}$$

In other words, calculating the cumulative hazard of the copula at $(p,q)$ is the same as calculating the cumulative hazard of the original distribution at the point $(F^{-1}(p), G^{-1}(q))$. We would like to estimate the cumulative hazard of the copula, and the equation above seems to suggest that we could do so by looking at the Nelson-Aalen estimator for the integrated hazard of the distribution $J$ evaluated at the appropriate quantiles of $F$ and $G$. If $F$ and $G$ are unknown, we would replace $(F^{-1}(p), G^{-1}(q))$ with the quantiles of the respective Kaplan-Meier estimates of the marginals $F$ and $G$ (See Section IV.3.1 in [1]).

Now, suppose for the moment that $F$ and $G$ are known. Define the pseudo-observations $(X_i^{\#}, Y_i^{\#}) := (F(X_i), G(Y_i))$, so that the distribution of

$(X_i^{\#}, Y_i^{\#})$ is the copula $C$. The observable region then becomes

$$\xi_i^{\#} = \{(F(x), G(y)) : (x, y) \in \xi_i\},$$

and the new filtration (now defined on $[0,1]^2$) is $\mathcal{F}^{\#}(p,q) = \mathcal{F}(F^{-1}(p), G^{-1}(q))$. Using Equation (3), we have

$$\hat{H}^{(n)c}(p,q) = \sum_{\{i:(X_i^{\#},Y_i^{\#})\le(p,q),(X_i^{\#},Y_i^{\#})\in\xi_i^{\#}\}} (Z_n^{\#}(X_i^{\#}, Y_i^{\#}))^{-1}. \tag{12}$$

But note that

$$\begin{aligned} Z_n^{\#}(u,v) &= \sum_{i=1}^{n} I_{\{(X_i^{\#},Y_i^{\#})\in E_{u,v}\}} I_{\{(u,v)\in\xi_i^{\#}\}} \\ &= \sum_{i=1}^{n} I_{\{(X_i,Y_i)\in E_{(F^{-1}(u),G^{-1}(v))}\}} I_{\{(F^{-1}(u),G^{-1}(v))\in\xi_i\}} \\ &= Z_n(F^{-1}(u), G^{-1}(v)), \end{aligned}$$

so

$$\begin{aligned} (12) &= \sum_{\{i:(X_i,Y_i)\le(F^{-1}(p),G^{-1}(q)),(X_i,Y_i)\in\xi_i\}} (Z_n(X_i,Y_i))^{-1} \\ &= \hat{H}^{(n)}(F^{-1}(p), G^{-1}(q)), \end{aligned} \tag{13}$$

and hence we have an empirical analogue of (11). More generally, for $D \subseteq [0,1]^2$, let $(F^{-1}, G^{-1})(D) := \{(F^{-1}(p), G^{-1}(q)) : p, q \in D\}$, in which case we can estimate $H_D^C$ with $\tilde{H}_D^{(n)C} = \hat{H}^{(n)}(F^{-1}, G^{-1})(D)$. The asymptotic normality of $\sqrt{n}(\tilde{H}_D^{(n)C} - H_D^C)$ follows from Remark 4.9.

If $F$ and $G$ are unknown, the next step is to estimate the quantiles $F^{-1}$ and $G^{-1}$. We will do so by taking the appropriate quantiles of the Kaplan-Meier estimators $\hat{F}$, $\hat{G}$ of $F$ and $G$ respectively; then our estimator will be $\hat{H}^{(n)C} = \hat{H}^{(n)}(\hat{F}^{-1}, \hat{G}^{-1})$, where $\hat{F}^{-1}(p) = \inf\{s : \hat{F}(s) \ge p\}$ and $\hat{G}^{-1}(q)$ is defined in a similar manner. It is important to note that $A_{(\hat{F}^{-1}(p),\hat{G}^{-1}(q))}$ is still a $*$-anti-cloud: indeed, since $\hat{F}(s)$ and $\hat{G}(t)$ are adapted to $\mathcal{F}^{(1)}(s)$ and $\mathcal{F}^{(2)}(t)$ respectively,

$$\begin{aligned} (s,t) \in A_{(\hat{F}^{-1}(p),\hat{G}^{-1}(q))} &= \left\{s \le \hat{F}^{-1}(p)\right\} \cap \left\{t \le \hat{G}^{-1}(q)\right\} \\ &= \left\{\hat{F}(u) < p,\ \forall\ u < s\right\} \cap \left\{\hat{G}(v) < q,\ \forall\ v < t\right\} \\ &\in \mathcal{F}_s^{(1)} \vee \mathcal{F}_t^{(2)} \subseteq \mathcal{F}_{(s,t)}^{*}. \end{aligned}$$

An immediate consequence of this is that our estimator still maintains a pseudo-strong martingale structure. Moreover, under Assumptions 4.2, 4.4 and 4.3, the composition map

$$\begin{aligned}(Z_n, {N^{(n)}}^{\xi}, Z_{n,1}, {N_1^{(n)}}^{\xi_1}, Z_{n,2}, {N_2^{(n)}}^{\xi_2}) &\xrightarrow{\phi_1} (Z_n, {N^{(n)}}^{\xi}, \hat{H}_1^{(n)}, \hat{H}_2^{(n)}) \\ &\xrightarrow{\phi_2} (Z_n, {N^{(n)}}^{\xi}, \hat{F}, \hat{G}) \\ &\xrightarrow{\phi_3} (\hat{H}^{(n)}, \hat{F}, \hat{G}) \\ &\xrightarrow{\phi_4} (\hat{H}^{(n)}, \hat{F}^{-1}, \hat{G}^{-1}) \\ &\xrightarrow{\phi_5} \hat{H}^{(n)} \circ (\hat{F}^{-1}, \hat{G}^{-1}) \qquad (14)\end{aligned}$$

is Hadamard-differentiable tangentially to $C[0,1]^2 \times BV_M[0,1]^2 \times (C[0,1] \times BV_M[0,1])^2$ on a domain analogous to map (5): the maps $\phi_1$ and $\phi_3$ are Hadamard-differentiable as a consequence of Theorem 4.5, $\phi_2$'s Hadamard differentiability comes from Lemma 3.9.30 in [13], and $\phi_4$ and $\phi_5$ are Hadamard-differentiable by Example 3.9.24 and Lemma 3.9.25, respectively, in [13] again. Hence, the process $\sqrt{n}(\hat{H}_{\cdot}^{(n)C} - H_{\cdot}^{C})$ converges in $l^\infty[0,\tau] = l^\infty([0,\tau_1] \times [0,\tau_2])$ to a mean zero Gaussian limit for every $(\tau_1, \tau_2)$ such that $S(F^{-1}(\tau_1), G^{-1}(\tau_2)) = \overline{C}(\tau_1, \tau_2) > 0$.

**An example: The FGM copula**

We can apply the preceding results to Farlie-Gumbel-Morgenstern (FGM) copulas. The FGM copula with parameter $\theta$ is defined by $C_\theta(u,v) = uv + \theta uv(1-u)(1-v)$. Suppose $C_{J_1}$ and $C_{J_2}$ are two FGM copulas with parameters $\theta_{J_1}$ and $\theta_{J_2}$, respectively, and we are interested in testing $H_\circ : C_{J_1} = C_{J_2}$ vs. $H_1 : C_{J_1} < C_{J_2}$. For FGM copulas, this is equivalent to testing $H'_\circ : \theta_{J_1} = \theta_{J_2}$ vs. $H'_1 : \theta_{J_1} < \theta_{J_2}$. This test can be performed by 'translating' these conditions to those of a test of hazard rate order: we look for a specific set $A$ under which $\theta_{J_1} < \theta_{J_2}$ is equivalent to $h^{C_{J_1}}(u,v) < h^{C_{J_2}}(u,v)$ for $(u,v) \in A$, and then we can apply the results of the preceding sections.

We have $\overline{C}_\theta(u,v) = (1-u)(1-v) + \theta uv(1-u)(1-v)$, as well as $\frac{\partial^2 C}{\partial u \partial v} = 1 + \theta(1-2u)(1-2v)$. Then $h^{C_{J_1}}(u,v) < h^{C_{J_2}}(u,v)$ becomes

$$\frac{1 + \theta_{J_1}(1-2u)(1-2v)}{(1-u)(1-v) + \theta_{J_1} uv(1-u)(1-v)} < \frac{1 + \theta_{J_2}(1-2u)(1-2v)}{(1-u)(1-v) + \theta_{J_2} uv(1-u)(1-v)},$$

which is equivalent to

$$(\theta_{J_2} - \theta_{J_1})uv < (\theta_{J_2} - \theta_{J_1})(1-2u)(1-2v)$$

after cross-multiplying and rearranging terms. Under the alternative hypothesis $\theta_{J_1} < \theta_{J_2}$, the previous inequality becomes

$$uv < (1-2u)(1-2v),$$

and finally

$$1 - 2u - 2v + 3uv > 0.$$

This means that, under $H_1' : \theta_{J_1} < \theta_{J_2}$, we have $h^{C_{J_2}}(u,v) - h^{C_{J_1}}(u,v) > 0$ for every $(u,v)$ lying in the set $\{(u,v) : 1-2u-2v+3uv > 0\}$. Now in order to ensure that $\overline{C_{J_i}}$ is uniformly bounded below, $i = 1,2$, define $A := [0,\tau] \cap \{(u,v) : 1-2u-2v+3uv > 0\}$ for any $\tau = (\tau_1, \tau_2)$ with $\tau_i < 1$, $i = 1,2$.

We consider two scenarios. First, suppose the distributions $J_1$ and $J_2$ have the same (known) marginal distributions $F$ and $G$. Since $\int_A (h^{C_{J_2}} - h^{C_{J_1}}) d\lambda = H_A^{C_{J_2}} - H_A^{C_{J_1}}$, an appropriate test statistic would be (letting $N = n + m$)

$$\sqrt{\frac{nm}{N}} \left( \tilde{H}_A^{(n)C_{J_2}} - \tilde{H}_A^{(n)C_{J_1}} \right),$$

which would get large if $\theta_{J_2} > \theta_{J_1}$.

More realistically, suppose the marginal distributions of $J_1$ and $J_2$ are unknown. If $\theta_{J_2} > \theta_{J_1}$, then $H_{A_u}^{C_{J_2}} - H_{A_u}^{C_{J_1}} > 0$ for every $u \in A$. In this case, an appropriate test statistic is $W_A^{(N)} = \sup_{u \in A} \sqrt{\frac{nm}{N}} (\hat{H}_{A_u}^{(m)C_{J_2}} - \hat{H}_{A_u}^{(n)C_{J_1}})$, which would get large if $\theta_{J_2} > \theta_{J_1}$. We need to consider two different cases for calculating the critical values for our test statistic.

If it is reasonable to assume that the marginals from the two distributions are equal, we proceed with the bootstrapping as indicated in Subsection 5.2. Under $H_\circ$, we know that $W_A^{(N)}$ converges to the sup over $A$ of a mean-zero Gaussian limit $U^{C_{J_1}}$, and our test of level $\alpha$ would reject $H_\circ$ whenever $W_A^{(N)} > \ddot{w}_N$, where $\ddot{w}_N = \inf \left\{ t : P(\ddot{W}_A^{(N)} > t) \leq \alpha \right\}$.

On the other hand, if the assumption of equal marginals is not justified, we can no longer pool the samples. In this case we will start by bootstrapping each sample separately. Let $\ddot{H}_{A_u}^{(n)C_{J_1}}$ and $\ddot{H}_{A_u}^{(m)C_{J_2}}$ denote the bootstrapped Nelson-Aalen estimators. Now let

$$\ddot{U}^{(n,m)} = \sqrt{\frac{nm}{N}} \left( \ddot{H}^{(m)C_{J_2}} - \ddot{H}^{(n)C_{J_1}} - (\hat{H}^{(m)C_{J_2}} - \hat{H}^{(n)C_{J_1}}) \right),$$

which can be rewritten as

$$\sqrt{\frac{n}{N}} \left( \sqrt{m} (\ddot{H}^{(m)C_{J_2}} - \hat{H}^{(m)C_{J_2}}) \right) - \sqrt{\frac{m}{N}} \left( \sqrt{n} (\ddot{H}^{(n)C_{J_1}} - \hat{H}^{(n)C_{J_1}}) \right).$$

Assuming that $\frac{n}{n+m} \to \lambda$ as $n, m \to \infty$, the last expression converges conditionally given the original samples to $\sqrt{\lambda} U^{C_{J_2}} - \sqrt{1-\lambda} U^{C_{J_1}}$ by Lemma 5.2, where $U^{C_{J_1}}$ and $U^{C_{J_2}}$ are the independent mean-zero Gaussian limits of $\sqrt{n}(\hat{H}^{(n)C_{J_1}} - H^{C_{J_1}})$ and $\sqrt{m}(\hat{H}^{(m)C_{J_2}} - H^{C_{J_2}})$ respectively. Under the null hypothesis, $W_A^{(N)}$ has the same limiting distribution as $\ddot{W}_A^{(n,m)} := \sup_{u \in A} \ddot{U}_{A_u}^{(n,m)}$. Then we have a test of asymptotic level $\alpha$ if we reject the null hypothesis whenever $W_A^{(N)} > \ddot{w}_{n,m}$, where $\ddot{w}_{n,m} = \inf \left\{ t : P(\ddot{W}_A^{(n,m)} > t) \leq \alpha \right\}$.

# References

[1] Andersen, P. K., Borgan, O., Gill, R. D. and Keiding, N., *Statistical Models Based on Counting Processes*, Springer Series in Statistics, Springer-Verlag, New York, 1993.

[2] Carabarin Aguirre, A., *Set-Indexed Survival Analysis with Generalized Censoring*, doctoral thesis, University of Ottawa, 2008.

[3] David-Zadeh, O., *Multi-Parameter Survival Analysis*, doctoral thesis, Bar-Ilan University, 2008.

[4] Hildebrandt, T. H., *Introduction to the Theory of Integration*, Academic Press Inc., New York, 1963.

[5] Hougaard, Philip, *Analysis of Multivariate Survival Data*, Springer, New York, 2000.

[6] Ivanoff, B. G. and Merzbach, E., *Set-Indexed Martingales*, Chapman and Hall/CRC Press, Boca Raton, 2000.

[7] Ivanoff, B. G. and Merzbach, E., Set-Indexed Markov Processes, *Canadian Mathematical Society Conference Proceedings* 26, 2000, 217-232.

[8] Ivanoff, B. G. and Merzbach, E., Random censoring in set-indexed survival analysis, *The Annals of Applied Probability* 12, 2002, 944-971.

[9] Ivanoff, B. G. and Merzbach, E., Random clouds and an application to censoring in survival analysis, *Stochastic Processes and their Applications* 111, 2004, 259-279.

[10] Kalbfleisch, J. D. and Prentice, R. L. *The Statistical Analysis of Failure Time Data*, Wiley, New Jersey, 2002.

[11] Pons, O., A test of independence between two censored survival times, *Scand. J. Statist.* 13, 173-185 (1986).

[12] Pons, O. and Turckheim, E. de, Tests of independence for bivariate censored data based on the empirical joint hazard function, *Scand. J. Statist.* 18, 21-37 (1989).

[13] van der Vaart, Aad W. and Wellner, Jon A., *Weak Convergence and Empirical Processes*, Springer series in statistics, 1996.